\documentclass{amsart}

\usepackage{amsmath}
\usepackage{amsfonts}
\usepackage{amssymb}
\usepackage{comment}
\usepackage{epsfig}
\usepackage{psfrag}
\usepackage{mathrsfs}
\usepackage{amscd}
\usepackage[all]{xy}
\usepackage{rotating}
\usepackage{lscape}
\usepackage{amsbsy}
\usepackage{verbatim}
\usepackage{moreverb}

\usepackage{color}

\newtheorem{theorem}{Theorem}[section]

\newtheorem{conj}[theorem]{Conjecture}

\newtheorem{definition}[theorem]{Definition}

\theoremstyle{remark}
\newtheorem{remark}[theorem]{Remark}
\newtheorem{example}[theorem]{Example}

\DeclareMathOperator{\Ext}{Ext}

\DeclareMathOperator{\Hom}{Hom}

\DeclareMathOperator{\Spec}{Spec}

\def\1bord{1\mathrm{Bord}}
\def\2bord{2\mathrm{Bord}}
\def\3bord{3\mathrm{Bord}}

\newcommand{\inv}{{}^{-1}}

\def\cA{\mathcal A}\def\cB{\mathcal B}\def\cC{\mathcal C}\def\cD{\mathcal D}
\def\cF{\mathcal F}\def\cH{\mathcal H}
\def\cK{\mathcal K}\def\cL{\mathcal L}
\def\cM{\mathcal M}\def\cN{\mathcal N}\def\cO{\mathcal O}

\def\cX{\mathcal X}

\def\AA{\mathbb A}\def\CC{\mathbb C}
\def\GG{\mathbb G}

\def\PP{\mathbb P}
\def\QQ{\mathbb Q}\def\RR{\mathbb R}

\def\ZZ{\mathbb Z}

\def\Gm{{\mathbb G}_m}

\def\QC{{\mathcal Q}{\mathcal C}}

\def\ot{\otimes}

\def\ot{\otimes}

\def\fM{\mathfrak M}

\def\wt{\widetilde}
\def\D{{\mathcal D}}

\newcommand{\on}{\operatorname}
\newcommand{\Vect}{\on{Vect}}
\newcommand{\Rep}{\on{Rep}}
\newcommand{\module}{-\on{mod}}

\newcommand{\Perf}{\on{Perf}}

\newcommand{\wh}{\widehat}

\def\oo{\infty}

\newcommand{\bs}{\backslash}

\newcommand\Cx{{\CC^\times}}

\def\Coh{\on{Coh}}

\def\Sym{\on{Sym}}

\def\End{\on{End}}

\def\fg{\mathfrak g}
\def\fb{\mathfrak b}

\def\ft{\mathfrak t}

\def\fM{\mathfrak M}

\def\St{St}

\def\fh{\mathfrak h}
\def\mf{\mathfrak}

\newcommand{\cSpec}{{\mathcal Spec}}

\newcommand{\Shv}{{\mathcal Shv}}

\newcommand{\Mob}{{\mathcal M\ddot{o}b}}
\newcommand{\Bun}{{\mathcal Bun}}
\newcommand{\Loc}{{\mathcal Loc}}
\newcommand{\Conn}{{\mathcal Conn}}
\newcommand{\Op}{{\mathcal Op}}
\newcommand{\Higgs}{{\mathcal Higgs}}
\newcommand{\Perv}{{\mathcal Perv}}

\newcommand{\fdet}{{\mf d}{\mf e}{\mf t}}

\begin{document}

\title{Betti Geometric Langlands}
\author{David Ben-Zvi and David Nadler}

\address{Department of Mathematics\\University of Texas\\Austin, TX 78712-0257}
\email{benzvi@math.utexas.edu}
 \address{Department of Mathematics\\University
  of California\\Berkeley, CA 94720-3840}
\email{nadler@math.berkeley.edu}

\maketitle

\newcommand{\GoG}{\displaystyle{\frac{G}{G}}}
\newcommand{\BoB}{\displaystyle{\frac{B}{B}}}
\newcommand{\BGB}{B\bs G/B}
\newcommand{\Gv}{{G^{\vee}}}
\newcommand{\Hv}{{H^{\vee}}}
\newcommand{\Bv}{{B^{\vee}}}

\newcommand{\Tv}{{T^{\vee}}}

\newcommand{\Pv}{{P^{\vee}}}

\newcommand{\Mv}{{M^{\vee}}}

\newcommand{\Badj}{{\mathbf{B^{\vee}}}}

\newcommand{\Gadj}{{\mathbf{G^{\vee}}}}

\newcommand{\cyl}{{Cyl}} 

\newcommand{\findim}{\mathit{fd}}





\begin{abstract}
We introduce and survey a Betti form of the geometric Langlands conjecture, parallel to the de Rham form developed by Beilinson-Drinfeld and Arinkin-Gaitsgory, and the Dolbeault form of Donagi-Pantev, and inspired by the work of Kapustin-Witten in supersymmetric gauge theory. The conjecture proposes an automorphic category associated to a compact Riemann surface $X$ and complex reductive group $G$ is equivalent to  a spectral category associated to the underlying topological surface $S$ and  Langlands dual group $\Gv$.
The automorphic category consists of suitable $\CC$-sheaves on the moduli  stack $\Bun_G(X)$ of $G$-bundles on $X$, 
 while the spectral category consists of  suitable $\cO$-modules on the character stack $\Loc_\Gv(S)$ of $\Gv$-local systems on $S$.
 The conjecture is compatible with and constrained by the natural symmetries of both sides
coming from modifications of bundles and local systems. 
 On the one hand, cuspidal Hecke eigensheaves in the de Rham and Betti sense are expected to coincide, so that one can view the Betti conjecture
 as  offering a different ``integration measure" on the same fundamental objects. On the other hand, the Betti spectral categories are more explicit than their de Rham counterparts and one might hope the conjecture is less challenging. The Betti program also enjoys symmetries coming from topological field theory: it is expected to extend to an equivalence of four-dimensional topological field theories, and in particular, the conjecture for closed surfaces is expected to reduce to the case of the thrice-punctured sphere. Finally, we also present ramified, quantum and integral variants of the conjecture,
and highlight connections to other topics, including representation theory of real reductive groups and quantum groups.
 \end{abstract}

\tableofcontents



\section{Introduction}


\subsection{Nonabelian Hodge theory}

Nonabelian Hodge theory on a smooth projective complex curve $X$, as formulated by Simpson (see e.g.~\cite{simpson}), 
studies three different moduli problems for  bundles for a complex reductive group $\Gv$:
\begin{enumerate}
\item[$\bullet$][deRham] $\Conn_\Gv(X)$: the moduli stack of flat $\Gv$-connnections on $X$,
\item[$\bullet$][Dolbeault] $\Higgs_\Gv(X)$: the moduli stack of $\Gv$-Higgs bundles on $X$,
\item[$\bullet$][Betti] $\Loc_\Gv(X)$: the moduli stack of $\Gv$-local systems on $X$.
\end{enumerate}
These carry a package of structures and relations generalizing the relations between de Rham, Dolbeault and Betti cohomology of smooth projective varieties including:
\begin{enumerate}
\item[$\bullet$] Given a real or integral form of $\Gv$, one can define a real or integral form of 
the Betti space.
\item[$\bullet$] The Betti space depends only on the homotopy type of $X$, and in particular, carries an action of the mapping class group of $X$.
\item[$\bullet$] The Riemann-Hilbert correspondence (from connections to their monodromy) provides an analytic identification from the de Rham space to Betti space. It induces an isomorphism of formal neighborhoods of points in the two spaces.
\item[$\bullet$] The de Rham space carries a nonabelian Hodge filtration, expressed via the Rees construction as a $\Gm$-equivariant family over $\AA^1$ with special fiber the Dolbeault space. 
\item[$\bullet$] The nonabelian Hodge theorem provides a trivialization of the Hodge filtration (a diffeomorphism between the Dolbeault and de Rham spaces) after passing from the stack to the corresponding moduli space of semistable objects. It induces an isomorphism of formal neighborhoods of points in the two spaces.
\item[$\bullet$] The de Rham space carries a flat (nonabelian Gauss-Manin)  connection over the moduli of curves. However, the connection does not integrate algebraically to a parallel transport:  de Rham spaces for distinct curves are not isomorphic.
\end{enumerate}

The Betti spaces are more elementary in construction than their de Rham and Dolbeault counterparts. Indeed, they are global quotients of  affine derived complete intersections. Their manifestly topological nature provides a ready source of examples of topological field theories. Their cohomology has been the subject of great recent interest~\cite{HRV,HLRV,dCHM}. They are central objects in the theory of cluster varieties~\cite{FG,FST,GS,STWZ}. 
 Our goal in this paper is to suggest a role for Betti spaces in the geometric Langlands program.

\subsection{Geometric Langlands Program}
The geometric Langlands program provides a nonabelian, global and categorical form of harmonic analysis.
We fix a complex reductive group $G$ and study the moduli stack $Bun_G(X)$ of $G$-bundles on $X$. This stack comes equipped with a large commutative symmetry algebra: for any point $x\in X$ we have a family of correspondences acting on $Bun_G(X)$ by modifying $G$-bundles at $x$. 
The goal of the geometric Langlands program is to simultaneously diagonalize the action of Hecke correspondences on suitable categories of sheaves on $Bun_G(X)$. One can ask to label the common eigensheaves (Hecke eigensheaves) by their eigenvalues (Langlands parameters), or more ambitiously, to construct a Fourier transform identifying categories of sheaves with dual categories of sheaves on the space of Langlands parameters. 

The kernels for Hecke modifications are bi-equivariant sheaves on the loop group $G(\cK)$,  $\cK= \CC((t))$,
with respect to the arc subgroup $G(\cO)$, $\cO =\CC[[t]]$.
The underlying double cosets 
are in bijection with irreducible representations of the Langlands dual group:
$$
\xymatrix{
G(\cO)\backslash G(\cK)/G(\cO)\ar@{<->}[r] &  \on{IrrRep}(\Gv)
}
$$
The geometric Satake theorem lifts this bijection to an equivalence of tensor categories 
$$\Perv(G(\cO)\backslash G(\cK)/G(\cO))\simeq \on{Rep}(\Gv)
$$ between equivariant perverse sheaves on the affine Grassmannian $Gr = G(\cK)/G(\cO)$ and finite dimensional representations of $\Gv$. This leads to the geometric notion of Langlands parameter: a $\Gv$-local system on $X$ provides a Hecke eigenvalue in that it defines 
a tensor functor $\on{Rep}(\Gv)\to \Vect$, for each point $x\in X$. 

\begin{conj}[Core Geometric Langlands]
For any irreducible $\Gv$-local system $E$ on $X$, there exists a perverse sheaf $\rho_E\in \Perv(Bun_G(X))$ 
with the structure of $E$-Hecke eigensheaf. 
\end{conj}

In order to lift this object-wise correspondence to a spectral decomposition, Beilinson and Drinfeld suggested we use the Riemann-Hilbert correspondence to replace perverse sheaves with the corresponding regular holonomic $\D$-modules. The category of all $\D$-modules provides a powerful geometric substitute for classical function spaces of harmonic analysis, in which we replace generalized functions by the system of all linear PDE with polynomial coefficients which they satisfy. More precisely, for a stack $M$ we let $\D(M)$ denote the dg-enhanced derived category of quasicoherent $\D$-modules on $X$, which is a differential graded category.\footnote{We view dg categories without further mention through the lens of homotopical algebra, i.e., as objects of the corresponding symmetric monoidal $\infty$-category.}
Roughly speaking, we then wish to identify $\D(Bun_G(X))$ with the dg category $\QC(\Conn_\Gv(X))$ of quasicoherent sheaves on the de Rham moduli space for the dual group. The singular nature of the de Rham space, or dually the ``noncompact" (specifically, non-quasicompact) nature of $Bun_G(X)$, forces us to be careful about regularity conditions on $\Conn_\Gv(X)$, or dually growth conditions on $Bun_G(X)$. These issues were solved by Arinkin-Gaitsgory~\cite{AG}, who introduced the dg category
$\QC^!_\cN(\Conn_\Gv(X))$ of ind-coherent sheaves with nilpotent singular support,
and proposed the following refined conjecture:

\begin{conj}[de Rham Geometric Langlands]
There is an equivalence of dg categories
$$\xymatrix{\D(Bun_G(X)) \simeq \QC^!_\cN(\Conn_\Gv(X))}$$
compatible with actions of Hecke functors. 
In particular, skyscrapers on $\Conn_\Gv(X)$ correspond to Hecke eigensheaves.
\end{conj}

We refer to~\cite{outline} for an overview of spectacular recent progress towards this conjecture.
The de Rham geometric Langlands program carries analogues of the familiar structures on the de Rham space $\Conn_\Gv(X)$.
The category $\D(Bun_G(X))$ has a flat connection over the moduli of curves (it forms a crystal of categories), and the conjecture is compatible with Gauss-Manin connections. Note that the connection is not integrable: the categories for distinct curves are never equivalent.

The category $\D(Bun_G(X))$ also carries a Hodge filtration: by degenerating differential operators to symbols, we have a $\Gm$-equivariant family of categories over $\AA^1$ with special fiber $\QC(T^*Bun_G(X))$, the category of quasi coherent {\em Higgs sheaves} on $Bun_G(X)$. The compatibility of the de Rham conjecture with the Hodge filtration (i.e., its semiclassical asymptotics) was studied by Arinkin~\cite{arinkin thesis}. A more subtle aspect of the Hodge filtration (identified in these terms in~\cite{simpson opers}) is the presence of a distinguished subvariety $\Op_\Gv(X)\subset\Conn_\Gv(X)$, the {\em opers} of Beilinson-Drinfeld, for which the geometric Langlands conjecture was proved in~\cite{BD}: the structure sheaf $\cO_{\Op_\Gv(X)}$ corresponds to the $\D$-module $\D$ itself (up to choices of spin structures). Finally, the associated graded (or ``classical limit") analog of the geometric Langlands conjecture is the following conjecture studied by Donagi and Pantev~\cite{DP Hitchin}, which we dub the Dolbeault form of Geometric Langlands:

\begin{conj}[Dolbeault Geometric Langlands]
There is an equivalence of dg categories
$$\xymatrix{\QC(T^*Bun_G(X)) \simeq \QC(T^*Bun_\Gv(X))}$$
\end{conj}

\begin{remark}
The conjecture is compatible with suitable Dolbeault Hecke functors. It yet remains to be suitably modified to account for singularities and non-compactness on both sides.
\end{remark}

The Dolbeault conjecture was proved in~\cite{DP Hitchin} over a dense open locus, where they reduce it to a form of the Fourier-Mukai transform for abelian varieties, applied to the fibers of Hitchin's integrable system. In particular, generic skyscrapers on $\QC(T^*Bun_\Gv(X))$ correspond to line bundles on smooth fibers of the Hitchin system. Following an idea of Donagi~\cite{donagi lecture}, Donagi and Pantev~\cite{DP Hodge} further explain how nonabelian Hodge theory on $Bun_G(X)$ (relating Higgs sheaves and $\D$-modules on $Bun_G(X)$) should directly relate the de Rham conjecture to the Dolbeault conjecture, a program they are pursuing with Simpson.



\subsection{Betti Geometric Langlands}
Our aim is to introduce a Betti form of the Geometric Langlands conjecture, in which the Betti space $\Loc_\Gv(X)$  takes the place of the de Rham space $\Conn_\Gv(X)$. Since closed points of the de Rham and Betti spaces are in canonical bijection and their formal neighborhoods are algebraically equivalent, the category generated by the core objects on the spectral side, the skyscrapers, will be the same in the Betti and de Rham versions.
On the automorphic side, we seek a Betti category containing the core automorphic objects, the perverse sheaves on $Bun_G(X)$ which are Hecke eigensheaves. As we discuss below, eigensheaves are expected to have characteristic variety contained in Laumon's global analogue of the nilpotent cone, namely the zero fiber $\cN_X=Hitch\inv(0)$ of the Hitchin fibration. (This is compatible with existing constructions in the de Rham setting,
as well as the motivating fact that eigensheaves in the Dolbeault setting are line bundles on Hitchin fibers, which all have 
 $\cN_X$ as the support of their conical limit.)
We thus propose the category  $\Shv_\cN(Bun_G(X))$ of {\em nilpotent sheaves}, complexes of sheaves of $\CC$-vector spaces, which are locally constant in codirections which are not nilpotent. Finite rank nilpotent sheaves are automatically constructible, however we require no finiteness conditions so land outside of the traditional realm of constructible sheaves.

\begin{conj}[Betti Geometric Langlands]\label{Betti conjecture}
There is an equivalence of dg categories
$$\xymatrix{\Shv_\cN(Bun_G(X)) \simeq \QC^!_\cN(\Loc_\Gv(X))}$$
compatible with actions of Hecke functors.
\end{conj}

Here are some appealing features of the Betti conjecture:

\begin{enumerate}
\item[$\bullet$] It contains and promotes the Core Geometric Langlands Conjecture to an alternative categorical equivalence.
\item[$\bullet$] It has a natural integral form for nilpotent sheaves of abelian groups 
and ind-coherent sheaves on the integral character stack.
\item[$\bullet$] It predicts that the category of nilpotent sheaves depends only on the topology of the curve,  which appears far from obvious from the definition.
\item[$\bullet$] The Betti spectral category can be glued together from a pants decomposition of the curve. There is a
parallel conjectural  gluing formula for nilpotent sheaves, reducing the conjecture to a small number of basic building blocks.
\item[$\bullet$] It has a natural ramified extension involving parabolic structures on both sides.
\item[$\bullet$] It has a natural quantum version, relating twisted sheaves to a deformation of the spectral category built out of the representations of the quantum group. 
\item[$\bullet$] It has a natural extension to unoriented surfaces, which
 for the M\"obius strip is directly related to Langlands-Vogan-Soergel duality for representations of real  groups.
\item[$\bullet$] The Betti categories in genus one (elliptic character sheaves) have close ties to the representation theory of $p$-adic groups and double affine Hecke algebras.

\end{enumerate}

An overarching theme is that the Betti conjecture fits into topological field theory (compared to the de Rham conjecture which fits into conformal field theory).
 The spectral categories carry many structures from 4-dimensional topological field theory, which correspond to structures (old and new, established and conjectural) on the automorphic side. Among these we discuss domain walls (functoriality), surface operators (ramification data), Wilson lines (modifications at points), Verlinde loops (modifications along loops), and local operators (measuring singular support).

\begin{remark}[Opers]
Opers play a central role in the de Rham setting absent in  the Betti setting. 
The variety $\Op_\Gv(X)\subset \Conn_\Gv(X)$ of opers is not algebraic in the Betti space, though individual opers can be considered. On the automorphic side, this corresponds to the fact that the $\D$-module $\D\in \D({\Bun_G(X)})$ itself, the ``canonical coisotropic brane" of~\cite{KW}, does not have a Betti counterpart. In this regard, the situation is not symmetric:  the structure sheaf $\cO\in \QC^!_\cN(\Loc_\Gv(X))$ on the spectral side is a reasonable Betti object, whose dual is a  ``nilpotent Whittaker sheaf" $\mathit{\cF_\cO}\in \Shv_\cN(Bun_G(X))$. 
\end{remark}

\begin{remark}[Nilpotent cones]\label{nilpotent discussion}
The roles of nilpotent cones on the automorphic and spectral sides, though suggestively parallel, are not known to us to be related: 
the former is controlled by $H^1$ of the curve with coefficients in a $G$-bundle,
while the latter by $H^2$ of the curve with coefficients in a $G^\vee$-local system.

On the spectral side, the nilpotent cone  controls the behavior of sheaves at singular points,
appearing for the same reasons as in the de Rham setting, though there one considers all $\D$-modules without any bound on their singular support.

On the automorphic side, the nilpotent cone controls the behavior of sheaves at ``infinity",
playing the role of a Lagrangian skeleton or superpotential familiar in symplectic geometry.
One can view nilpotent sheaves as a form of ``partially wrapped" Lagrangian branes for the Hitchin system.
While it is beyond current technology to define the Fukaya category for $T^*Bun_G(X)$, we believe the category $\Shv_\cN(Bun_G(X))$ to be a good algebraic model for the category of $A$-branes from~\cite{KW}. 


In~\cite{KW}, Kapustin and Witten studied a topologically twisted form of maximally supersymmetric ($\cN=4$) super Yang-Mills theory in four dimensions, and related it to the geometric Langlands correspondence. 
The automorphic category proposed in~\cite{KW} is the category of D-branes in the topological A-model with target the Hitchin space $(\fM^H_G(X),\omega_K)\supset (T^*\fM_G^s, \omega)$, the moduli space of semistable $G$-Higgs bundles on $X$, with its symplectic structure which restricts to the standard form on the cotangent to the moduli space of stable $G$-bundles.

A beautiful feature of the Kapustin-Witten A-model picture of automorphic sheaves is that there is a large supply of obvious automorphic objects (Hecke eigensheaves). Namely, any rank one local system on a Lagrangian torus, given by a smooth fiber of $Hitch$, defines such an object. Under mirror symmetry (T-duality along the Hitchin fibration), this object is sent to a B-brane on the moduli space of $\Gv$-local systems on $X$, which is a skyscraper at a smooth point. One then finds~\cite{KW,Witten revisit} that the original A-brane is an eigenbrane for the 't Hooft operators, with eigenvalue given by the corresponding $\Gv$-local system.

Kapustin-Witten explain a connection between the A-model on $\fM^H_G$ and $\D$-modules on $Bun_G(X)$. The mechanism they propose involves constructing a $\D$-module out of an A-brane $F$ by considering homomorphisms $Hom(\cB_{cc},F)$ to $F$ out of a distinguished space-filling A-brane, the canonical coisotropic brane $\cB_{cc}$. These homomorphisms form a module for the endomorphisms of $\cB_{cc}$, which they identify in terms of differential operators on $Bun_G$. However the full ring of endomorphisms appears to be a ring of infinite-order differential operators, which are {\em entire} rather than polynomial functions of momenta. Dually, the B-model of the noncompact analytic space $\Loc_\Gv(X)_{an}$ has as full ring of endomorphisms of the structure sheaf the ring of analytic functions on the Betti, or equivalently de Rham, analytic space. Thus one needs to impose suitable growth conditions to obtain an equivalence of algebraic categories such as Betti or de Rham. One can view this as the role of the nilpotent cone on the automorphic side.

\end{remark}

In this article, we survey the developing ideas contributing to the Betti geometric Langlands conjecture. We will be informal, in particular suppressing all the prevalent $\infty$-categorical considerations. Our perspective is deeply influenced by the structure of topological field theory, primarily through the work of Kapustin-Witten~\cite{KW} and Lurie~\cite{jacob TFT}. Finally, we have attempted to include helpful references but have certainly missed many relevant works.

\subsection{Acknowledgements}

We would like to thank Adrien Brochier, David Jordan, Penghui Li, Toly Preygel and Zhiwei Yun for their collaborations on different aspects of this project.
We would like to acknowledge the National Science Foundation for its support, both through FRG grant ``In and Around Theory $\mathfrak X$" (DMS-1160461) and individual grants DMS-1103525 (DBZ) and DMS-1502178 (DN). We would also like to acknowledge that part of the work was carried out at MSRI as part of the program on Geometric Representation Theory.


\section{Two toy models}

\subsection{Mellin transform}
To illustrate the relation between the de Rham and Betti conjectures, let us consider the Fourier transform on multiplicative groups (Mellin transforms).

Recall the classical Mellin transform:  
given a function $f(z):\RR_{>0} \to \RR$, its Mellin transform
$$
\xymatrix{
\displaystyle\wh{f}(s)=\int f(z) z^s \frac{dz}{z} :\RR\ar[r] & \RR
}$$ 
provides the coefficients of the expansion of $f(z)$ in terms of the characters 
$$
\xymatrix{
z^s=e^{sx} :\RR_{>0} \ar[r] &  \RR & s\in \RR
}$$ 
The characters $z^s = e^{sx}$ are homogeneous of degree $s$, and form a ``basis" of eigenfunctions of the invariant differential operator $\partial=z\frac{d}{dz}$. The Mellin transform ``diagonalizes" the operator $\partial$, while exchanging the multiplication by $z$ with translation in $s$.

Now let us consider an algebraic realization of the Mellin transform in several variables. 

Let $\Lambda$ be a lattice, with dual lattice $\Lambda^\vee$, and consider the torus $T=\Cx \otimes_\ZZ \Lambda$
with coweight lattice $\Lambda = \Hom(\Cx, T)$,
weight lattice $\Lambda^\vee = \Hom(T, \Cx)$, Lie algebra $\ft=\CC\otimes_\ZZ \Lambda$
and dual $\ft^*=\CC\otimes_\ZZ \Lambda^\vee$.

On the one hand, 
let $z_i\in \Lambda^\vee$ be a basis, and
consider the regular functions
 $\CC[T]= \CC[z_i^{\pm 1}]$ and
the algebra of differential operators $$\cD_T=\CC[z_i^{\pm 1}]\langle \partial_i\rangle/ \{ \partial_i z_j=  z_i (\delta_{ij}+ \partial_i)\}$$ 
where $\partial_i=z_i \frac{\partial}{\partial z_i} \in\Lambda \subset \ft$ are a basis of $T$-invariant vector fields. 

On the other hand, write $\xi_i  \in \Lambda \subset \ft$ for a basis,
and consider the regular functions
 $\CC[\ft^*]= \CC[\xi_i]$ and
the algebra of finite difference operators 
$$ \Delta_{\ft^*}= \CC[\xi_i]\langle \sigma_i^{\pm 1}\rangle / \{ \xi_i\sigma_j=\sigma_j (\delta_{ij}+\xi_i) \}
$$
where  $\sigma_i\in \Lambda^\vee$ are a basis of the shift operators.

There is an evident isomorphism 
$$
\xymatrix{
\cD_T\simeq \Delta_{\ft^*}
& z_i \longleftrightarrow \sigma_i
& \partial_i \longleftrightarrow \xi_i  
}$$
and thus an equivalence of abelian categories
$$
\D_T\module \simeq  \QC(\ft^*)^{\Lambda^\vee}
$$ 
between $\D_T$-modules and $\Lambda^\vee$-equivariant quasicoherent sheaves on $\ft^*$.
The equivalence matches the eigensystems of $\partial_i$ in the form of the irreducible local systems 
$$
\xymatrix{
\cL_\lambda = \cD_T/\D_T \langle \partial_i - \lambda_i \rangle = \{ \partial_i f= \lambda_i f\}
&
 \lambda \in \ft^*
 }$$
with the difference modules 
$$
\xymatrix{
\bigoplus_{\mu\in \Lambda^\vee} \cO_{\lambda + \mu}
&
 \lambda \in \ft^*
 }$$ 
given by the $\Lambda^\vee$-equivariantizations of the skyscrapers $\cO_\lambda$,
in particular expressing the gauge equivalence of the eigensystems $\cL_\lambda$ and $\cL_{\lambda+\mu}$, for $\mu\in \Lambda^\vee$.

The equivalence realizes arbitrary $\D_T$-modules as ``direct integrals" of the eigensystems $\cL_\lambda$
in the sense that $\Lambda^\vee$-equivariant quasicoherent sheaves on $\ft^*$ are ``direct integrals" of the skyscrapers
$\cO_\lambda$.
Passing to the 
locally defined solutions $Cz^{\lambda}$ of
the eigensystems $\cL_\lambda$, we arrive back at the  Mellin transform in several variables.

Now instead of $\D_T$-modules,  let us consider the abelian category $\Loc(T)$ of local systems on $T$ of arbitrary rank
with its evident equivalent descriptions
  $$
  \xymatrix{
  \Loc(T)=\CC[\pi_1(T)]\module= \CC[\Lambda]\module \simeq \QC(T^\vee) 
  }
  $$ 
  where we introduce  the dual torus
  $T^\vee=\Cx \otimes_\ZZ \Lambda^\vee$
with coweight lattice $\Lambda^\vee = \Hom(\Cx, T)$ and
weight lattice $\Lambda = \Hom(T, \Cx)$.

The irreducible local system $L_\sigma$ of monodromy $\sigma\in T^\vee$ corresponds to the skyscraper $\cO_\sigma$ at the point $\sigma \in T^\vee$.
Arbitrary local systems are ``direct integrals" of the irreducibles $\cL_\sigma$
in the sense that quasicoherent sheaves on $T^\vee$ are ``direct integrals" of the skyscrapers
$\cO_\sigma$.

The Riemann-Hilbert correspondence sends the $\D_T$-module $\cL_\lambda$ to
the local system $L_{\sigma}$ where $\sigma=\exp(\lambda)$. It extends to an equivalence of analytic moduli spaces 
$$
\xymatrix{
\exp:(\ft^*/\Lambda^\vee)_{an}\ar[r]^-\sim &  (T^\vee)_{an}
}
$$ 

While the ``small" irreducible objects $\cL_\lambda\in \D_T\module $ and $L_\sigma\in \Loc(T)$, where
$\sigma=\exp(\lambda)$,  are in natural correspondence,
  ``large" objects of the respective categories are of very different natures. For example, 
 consider the identity $e\in  T$,
 and the respective skyscrapers based there.
 On the one hand, we obtain  the $\D_T$-module of delta-functions
   $$
   \D_T/\D_T (z_i-1)\simeq \CC[\ft^*] \in \D_T\module
   $$ corresponding to the structure sheaf of $\ft^*$ with its standard $\Lambda^\vee$-equivariant structure. 
   On the other hand, we obtain the universal  local system represented by
   the regular $\CC[\pi_1(T)]$-module 
   $$ \CC[\pi_{1}(T)]  \simeq \CC[T^\vee] \in \Loc(T)
   $$ 
   corresponding to the structure sheaf of $T^\vee$. The endomorphisms of the first are scalars $\CC$,
   while the endomorphisms of the second are regular functions $\CC[T^\vee]$.


\subsection{Geometric class field theory on an elliptic curve}\label{elliptic toy}
Here we outline the basic shape of the Dolbeault, de Rham,  and Betti geometric Langlands equivalences in the case of line bundles $G=\Gv=GL_1$ on an elliptic curve $X=(E,0)$, as well as their quantum counterparts.

\subsubsection{Dolbeault}
Let $Jac(E, 0)=Pic^0(E, 0)\simeq E$ denote the moduli of line bundles of degree zero trivialized at $0\in E$. The Dolbeault space is $T^*Jac(E, 0)\simeq E\times \CC$ with its Hitchin integrable system given by projection to the second factor. The self-duality of the Jacobian induces a Fourier-Mukai auto-equivalence of the dg category 
$\QC(Jac(E, 0))$ of quasicoherent sheaves, 
and hence a fiberwise auto-equivalence of
$\QC(T^*Jac(E, 0))$. It exchanges a skyscraper on a fiber with a degree zero line bundle on the same fiber. 

\subsubsection{de Rham}
Let $Conn_{GL_1}(E,0)$ denote the moduli space of flat line bundles on $E$ trivialized at $0\in E$. 
Forgetting the connection realizes $Conn_{GL_1}(E,0)$ as an $\AA^1$-bundle over $Jac(E, 0)=Pic^0(E, 0)\simeq E$, specifically
  the unique nontrivial $\AA^1$-bundle, called
  the Serre surface. 
  The Fourier-Mukai transform for $\D$-modules of Laumon and Rothstein provides an equivalence
  of dg categories
$$
\D(Pic^0(E, 0))\simeq \QC(Conn_{GL_1}(E,0))$$
 It matches a flat line bundle $\cL \in \D_{Pic^0(E, 0)}\module$ with the corresponding skyscraper $\cO_{\cL}\in \QC(Conn_{GL_1}(E,0))$.

\subsubsection{Betti}
 Let $Loc_{GL_1}(E,0)$ denote the moduli of rank one local systems on $E$ trivialized at $0$. Taking monodromy
 around $\pi_1(Jac (E, 0))\simeq \ZZ\oplus \ZZ$ provides an equivalence 
 $$Loc_{GL_1}(E,0)\simeq \Spec \CC[\pi_1(Jac (E, 0))] \simeq \Cx\times \Cx 
 $$
so that we have an equivalence of  abelian categories
$$ Loc(Jac (E, 0))\simeq \QC(Loc_{GL_1}(E,0))$$ 
Note this is compatible with the calculations of wrapped microlocal sheaves and wrapped Lagrangian branes
on the Dolbeault space $T^*Jac(E, 0)$~\cite{wrapped loops, Nwms}.

\subsubsection{Quantum de Rham}
The $\AA^1$-bundle $Conn_{GL_1}(E,0) \to Jac(E, 0) \simeq E$ is the twisted cotangent bundle of $E$ associated to the divisor $0\in E$, i.e., sections are 1-forms on $E$ with simple poles at $0$ of residue 1. 
It carries a canonical algebraic symplectic form 
which can be quantized by twisted $\D$-modules on $E$.
For $1/\hbar=k\in \CC$, we have a category $\D_k\module$, where $$\D_k=\D_E(\cO(0)^k)$$ is the sheaf of differential operators on $E$ twisted by the $k$th power of the line bundle $\cO(0)$. This is a flat deformation as $\cO_E$-algebra of $\cO_{E^\flat}=\D_{\hbar=0}$. This deformation quantization is Fourier dual to the deformation of $\D_E$ to the sheaf of twisted differential operators on $E$:
the Fourier-Mukai transform for $\D$-modules has a twisted variant, due to Polishchuk-Rothstein, giving a derived equivalence
$$D(\D_{Jac E,k}\module)\simeq D(\D_{Jac E,1/k}\module)$$ degenerating to the above equivalence as $k\to 0$.

\subsubsection{Quantum Betti}\label{toy quantum Betti}
The variety $\Gm\times\Gm$ carries the algebraic symplectic form $\frac{dz}{z}\wedge \frac{dw}{w}$, which is analytically equivalent to the form on $E^\flat$. The corresponding Poisson bracket on functions is the classical limit $q\to 1$ of the q-Weyl algebra of difference operators (or quantum differential operators) on $\Gm$, $$\D_q(\Gm)=\ZZ\langle x,x\inv, \sigma,\sigma\inv\rangle/\sigma x=q x\sigma.$$ To see the Fourier dual of this quantization, let 
$\cL^{\log q}$ denote the $\Gm$-gerbe on $E$ associated to the $\log q$th power of the $\Gm$-bundle $\cL^\times$ of the line bundle $\cO(0)$, with class $[\cL^{\log q}]=q\in H^2(E,\Cx)\simeq \Cx$. Concretely, local systems on $E$ twisted by the gerbe $\cL^{\log q}$ are by definition local systems on the total space of $\cL^\times$ with fiber wise monodromy $q\in\Cx$. 
Such twisted local systems are ``level q" representations of a $\Cx$-central extension of $\pi_1(E)$, which we can describe explicitly by trivializing $\cL$ on $E\setminus 0$. The result is an abelian equivalence
$$\Loc_q(E) \simeq \D_q\module.$$


\section{Spectral side}

\subsection{Character stacks}
Given a topological surface $S$ and complex reductive group $\Gv$, 
we would like to study the  character stack
 of $G^\vee$-local systems 
$$Loc_\Gv(S)=[S,BG^\vee]$$
which we understand as a derived stack. 
Note this definition makes sense for $S$ an arbitrary homotopy type,
and many of our constructions will extend to this generality.

The character stack carries a natural action of the homotopy type $\mathit{Diff}(S)$. 
For $S$ a smooth oriented surface of genus $g>1$,
the action of  $\mathit{Diff}(S)$ factors through 
the mapping class group $MCG(S)=\pi_0(\mathit{Diff}(S))$.
More generally, a cobordism $M:S\leadsto S'$ of surfaces, or in fact any cospan $S\rightarrow M\leftarrow S'$ of spaces, 
defines a correspondence of character stacks
$$\xymatrix{\Loc_\Gv(S)& \ar[l] \Loc_\Gv(M) \ar[r]&\Loc_\Gv(S')}$$

\begin{remark}\label{integral form}
Reductive groups admit natural split forms over $\ZZ$ and  
as a result $\Loc_{G^\vee}(S)$ does as well.
All of the symmetries of  $\Loc_{G^\vee}(S)$ defined by symmetries of $S$ are defined over $\ZZ$.
\end{remark}

If we fix a point $s\in S$, the character stack admits an explicit description.
First, we have the smooth affine variety of representations
$$
\xymatrix{
\Rep_\Gv(S\setminus \{s\})=\Hom(\pi_1(S\setminus \{s\}), \Gv) 
}$$
Then the character stack has a global complete intersection 
presentation by group-valued Hamiltonian reduction
$$Loc_\Gv(S)\simeq(\Rep_\Gv(S\setminus \{s\}) \times_{\Gv} \{e\}))/\Gv$$ 
Thus one starts with $\Gv$-local systems on the punctured surface $S\setminus\{s\}$ trivialized at a base point, 
imposes that the monodromy around $s$ is the identity $e\in \Gv$, 
and then quotients by the adjoint action of $\Gv$ to forget the trivialization. Equivalently, one can first quotient by 
 the adjoint action of $\Gv$ to forget the trivialization, and then impose that the monodromy around $s$ is 
 conjugate to the identity $e\in \Gv$.
Thus one takes $\Gv$-local systems on the punctured surface $S\setminus\{s\}$ and imposes that their restriction
to the disk $D^2 \subset S$ around the puncture extends across the puncture
$$Loc_\Gv(S)\simeq\Loc_\Gv(S\setminus \{s\}) \times_{\Loc_{\Gv}(D^2\setminus x)} \Loc_{\Gv}(D^2)$$
In particular, for $S$ an oriented surface of genus $g$,
 we find the derived fiber product 
$$Loc_\Gv(S)\simeq  (\Gv)^{\times 2g}/\Gv  \times_{ \Gv/\Gv} \{e\}/\Gv $$
%
with $(\Gv)^{\times 2g}/\Gv \to \Gv/\Gv$ the usual product of commutators of monodromies,
and $\{e\}/\Gv \to \Gv/\Gv$ induced by the inclusion of the identity $e\in \Gv$. 

\begin{remark}
It is important to recognize that the character stack has a nontrivial derived structure, coming from the 
overdetermined nature of the group-valued moment map  $\Rep_\Gv(S\setminus \{s\}) \to G^\vee$.
At the linear level, this is recorded by the obstructions $H^2(S, ad(L))$, for a $G^\vee$-local system $L$.
\end{remark}

\subsubsection{Parabolic bundles}
Suppose the surface $S$ has boundary, fix a Borel subgroup $\Bv\subset \Gv$,
and let $H^\vee = B^\vee/[B^\vee, B^\vee]$ be the universal Cartan group.

We would also like to study the parabolic character stack  of maps of pairs
$$Loc_\Gv(S, \partial S)=[(S, \partial S),(BG^\vee, BB^\vee)]$$
This is an $\cX$-space in the sense of~\cite{FG}.

Going further, we have the class of the boundary monodromies
$$
\xymatrix{
Loc_\Gv(S, \partial S)\ar[r]  & \Loc_{H^\vee}(\partial S)
}
$$
and would often like to focus on  the unipotent fiber 
$$Loc^{u}_\Gv(S, \partial S) = [(S, \partial S),(BG^\vee, \tilde \cN^\vee/G^\vee)] \simeq Loc_\Gv(S, \partial S) \times_{ \Loc_{H^\vee}(\partial S)} (\{e\}/H^\vee)^{\pi_0(\partial S)}
$$
where we trivialize the classes of the boundary monodromies.
If we were to further trivialize the underlying boundary $H^\vee$-bundle by base-changing along
$\{e\}\to (\{e\}/H^\vee)^{\pi_0(\partial S)}$,
and twist by a central involution (see Remark~\ref{twisted}),
this would be an $\cA$-space in the sense of~\cite{FG}.

Finally, we  also have the completion of the unipotent fiber inside the parabolic character stack 
$$
\xymatrix{
Loc^{\hat u}_\Gv(S, \partial S)  
= [(S, \partial S),(BG^\vee, \widehat{\tilde \cN^\vee}/G^\vee)] \simeq Loc_\Gv(S, \partial S) \times_{ \Loc_{H^\vee}(\partial S)} 
(\widehat {\{e\}}/H^\vee)^{\pi_0(\partial S)}
}
$$

\subsubsection{Twisted groups}\label{twisted groups}
In the Langlands program, one is naturally led to local systems for non-constant group schemes.
In our setting, the definition and basic   properties of character stacks 
carry over easily to possibly non-trivial groups over $S$.

A natural example arises when $S$ is an unoriented surface.
Given an extension of groups 
$$\xymatrix{
1\ar[r] &  \Gv\ar[r] & G^L\ar[r] & \ZZ/2\ar[r]&  1
}$$
introduce the preimage $G^L_{-1} \subset G^L$  of the nontrivial element $-1\in \ZZ/2$.
 Define
 the corresponding twisted character stack
 $$\Loc_{G^L,or}(S) = \Loc_{G^L}(S) \times_{ \Loc_{\ZZ/2}(S)} \{\tilde S\}
 $$ 
parameterizing  $G^L$-local systems on $S$ equipped with an isomorphism between
 their induced $\ZZ/2$-local system and the orientation double cover $\tilde {S}\to S$.


\subsection{Examples}\label{examples}

\begin{example}[abelian case]
For $\Gv=\Tv$ a torus, the character stack is a product 
$$
\Loc_\Tv(S)\simeq B\Tv\times (\Tv\ot_\ZZ H^1(S,\ZZ)) \times {\mathfrak t}^{\vee}[-1]
$$
of the underlying classical character stack 
$$
\Loc_\Tv(S)_{cl} \simeq B\Tv\times (\Tv\ot_\ZZ H^1(S,\ZZ))
$$
which itself is a product of the classical character variety  $(\Tv\ot_\ZZ H^1(S,\ZZ))$
and the classifying stack $B\Tv$, along with the affine derived scheme
$${\mathfrak t}^{\vee}[-1]\simeq \{e\}\times_{\Tv} \{e\} = \Spec \Sym (\ft[1])$$
using the identification $ \ft\simeq (\ft^\vee)^*$.
\end{example}

\begin{example}[2-sphere]
For $S=S^2\simeq D^2 \amalg_{S^1} D^2$, we find
$$\Loc_\Gv(S^2)\simeq  \{e\}/\Gv\times_{\Gv/\Gv} \{e\}/\Gv\simeq \fg^{\vee}[-1]/\Gv$$ the derived self-intersection of the identity in the adjoint quotient.\end{example}

\begin{example}[3-sphere]
Similarly, for $S^3\simeq D^3 \amalg_{S^2} D^3$, we find
$$\Loc_\Gv(S^3)\simeq \{0\}/\Gv\times_{\fg^{\vee}[-1]/\Gv} \{0\}/\Gv\simeq \fg^{\vee}[-2]/\Gv$$  
the derived self-intersection of the origin within $\Loc_\Gv(S^2)$.
\end{example}

\begin{example}[cylinder] \label{cylinder example}
In the case of a cylinder $\cyl = S^1 \times [0,1]$ with boundary $\partial \cyl = S^1 \times \{0,1\}$,
we obtain the Grothendieck-Steinberg stack
$$
\xymatrix{
\Loc_\Gv(\cyl, \partial \cyl)\simeq \St_\Gv =   \Bv/\Bv \times_{\Gv/\Gv} \Bv/\Bv 
}
$$
$$
\xymatrix{
\simeq \{g\in \Gv, B_1,B_2\in \Gv/\Bv\; :\; g\in B_1\cap B_2\}/\Gv
}
$$ 
The unipotent version is the more familiar Steinberg stack
$$
\xymatrix{
\Loc_\Gv^u(\cyl, \partial \cyl)\simeq \St^u_\Gv =   \wt{\cN}^\vee/\Gv \times_{\Gv/\Gv} \wt{\cN}^\vee/\Gv}
$$
$$
\xymatrix{
\simeq \{g\in \Gv, B_1,B_2\in \Gv/\Bv\; :\; g\in B_1\cap B_2 , [g]_1 = [g]_2 = e\}/\Gv
}
$$
There is a nontrivial homotopical $S^1$-action  on $\St_\Gv$ given by rotating the cylinder.
\end{example}

\begin{example}[torus]\label{torus} For $S=T^2$, we find the derived commuting stack 
$$\Loc_\Gv(T^2)=\{g,h\in\Gv\; : \: gh=hg \}/\Gv$$ 
It carries an action of 
$$\mathit{Diff}(T^2)\simeq T^2 \rtimes SL_2(\ZZ)$$
where  the mapping class group $SL_2(\ZZ)$ permutes products of powers of $g$ and $h$, 
and the connected component $T^2$ acts by translations.
\end{example}

\begin{example}[pair of pants]
For $S$ the complement of two disjoint open disks in a closed disk, 
the fundamental group is free on two generators, and thus  we find
$$\Loc_\Gv(S,\partial S)=\{g,h\in \Gv,\; B_1,B_2,B_3\in \Gv/\Bv \; : \; g\in B_1,\;h\in B_2,\; gh\in B_3\}/\Gv.$$

\end{example}

\begin{example}[M\"obius strip] \label{mobius}

For an unoriented example, consider the M\"obius strip $\Mob$ where we quotient the cylinder $\cyl$ by the antipodal map. Equivalently, we also obtain $\Mob$ by removing an open disk from the projective plane  $\RR\PP^2$.

For simplicity, let us take the semi-direct product $G^L = G^\vee\rtimes \ZZ/2 $
where $\ZZ/2$ acts by an algebraic involution $\theta^\vee$.
In this case, the character stack
$$\Loc_{G^L,or}(\Mob,\partial \Mob)\simeq (\tilde \Gv \times_{\Gv} G^L_{-1})/G^\vee\simeq  \{g\in G^L_{-1}, B\in \Gv/\Bv\; : \; g^2\in B\}/\Gv$$
 is the Langlands parameter space ${\mathcal La}_\Gv^{\theta^\vee}$ of~\cite{reps}.
 The unipotent version $\Loc_{G^L,or}^u(\Mob, \partial \Mob)$ recovers the unipotent version  ${\mathcal La}_\Gv^{u,\theta^\vee}$
in which we impose that the group element $g$ is unipotent. 
It carries a homotopical $S^1$-action since rotating the cylinder commutes with the antipodal map.

\end{example}


\subsection{Spectral categories}
Let us continue with a topological surface $S$, complex reductive group $\Gv$, 
and  the resulting  character stack $\Loc_\Gv(S)$.

For the spectral side of the Betti Langlands conjecture, we would like to take an appropriate category
of $\cO$-modules on  $\Loc_\Gv(S)$.
 Due in particular to the presence of reducible local systems, $\Loc_\Gv(S)$ is singular, and so the natural candidate
 small dg categories
$\Coh(\Loc_\Gv(S))$ and $\Perf(\Loc_\Gv(S))$ of coherent sheaves and perfect complexes do not coincide.
Equivalently, 
 the natural  candidate cocomplete
 dg categories
$\QC^!(\Loc_\Gv(S))$ and $\QC(\Loc_\Gv(S))$ of ind-coherent sheaves and quasicoherent sheaves
do not coincide.

\begin{remark}[small and cocomplete categories]
One can go back and forth between the pair of small dg categories
$\Coh(\Loc_\Gv(S))$ and $\Perf(\Loc_\Gv(S))$  and the pair
of comcomplete dg categories
$\QC^!(\Loc_\Gv(S))$ and $\QC(\Loc_\Gv(S))$   by  taking ind-completions or compact objects.
More precisely, taking ind-categories defines a symmetric monoidal, colimit preserving equivalence between the $\infty$-category of small idempotent-complete stable $\infty$-categories with exact functors and the $\infty$-category of compactly-generated stable presentable $\infty$-categories with functors preserving colimits and compact objects.
We use this  to move freely between the two settings, though one must be careful when considering functors between 
cocomplete categories that do not preserve compact objects. 
\end{remark}

\begin{example}[abelian case]
Recall, for $\Gv=\Tv$ a torus, the character stack is a product 
$$
\Loc_\Tv(S)\simeq B\Tv\times (\Tv\ot_\ZZ H^1(S,\ZZ)) \times {\mathfrak t}^{\vee}[-1]
$$
with derived structure coming from the affine derived scheme
$${\mathfrak t}^{\vee}[-1]\simeq \{e\}\times_{\Tv} \{e\} = \Spec \Sym (\ft[1])$$

One can measure the difference between $\Coh(\Loc_\Tv(S))$ and $\Perf(\Loc_\Tv(S))$ by the action of the exterior algebra 
$\Lambda =  \Sym (\ft[1])$.
The singular support $ss(M)\subset \ft^\vee$ of a $ \Lambda$-module $M$ is the support of the  Koszul dual
graded module $M^\vee=\Ext_\Lambda(\CC_0,M)$ for the graded symmetric algebra  
$$S=\Ext_{\Lambda}(\CC_0,\CC_0)\simeq \Sym (\ft^\vee[-2])
$$ 
The singular support of $\Lambda$ itself, and hence of any nontrivial perfect $\Lambda$-module, is the origin $0\in \ft$. Thus nontrivial perfect $\Lambda$-modules are precisely the finitely-generated modules with  singular support
the origin $0\in \ft$.
\end{example}

The general constructions of Arinkin-Gaitsgory~\cite{AG} 
applied to the global complete intersection presentation of $ \Loc_\Gv(S)$
provide an action of the algebra  
$$
A= \Sym (\fg^{\vee}[-2])^\Gv\simeq \Sym (\fh^\vee[-2])^W
$$
by endomorphisms of the identity functor of  $\QC^!(\Loc_{\Gv}(S)).
$
Thus the endomorphisms of any object naturally form an $A$-algebra,
providing a notion of the $A$-support of the object
$$\xymatrix{
Supp_A(\cF) = Supp_{A}(End(\cF))\subset \fg^{\vee}//\Gv\simeq \fh^\vee//W
&
\cF\in \QC^!(\Loc_\Gv(S))
}
$$
In Section~\ref{local ops} below, we explain this structure in terms of local operators in topological field theory.

\begin{definition} 
1) An object $\cF\in \QC^!(\Loc_\Gv(S))$ has 
nilpotent  singular support if
 its  $A$-support 
 $$
 Supp_A(\cF) \subset \fg^\vee//G\simeq \fh^\vee//W
 $$  is  either the origin
  or empty.

2) Define the  small Betti spectral category 
 $$\Coh_\cN(\Loc_\Gv(S))\subset \Coh(\Loc_\Gv(S))$$
 and the large Betti spectral category
 $$\QC^!_\cN(\Loc_\Gv(S))\subset \QC^!(\Loc_\Gv(S))$$
  to be the respective full dg subcategories of objects   with nilpotent  singular support.
  
\end{definition}

\begin{remark}
Nilpotent  singular support includes the trivial singular support of quasicoherent complexes, and thus we also have 
$\Perf(\Loc_\Gv(S))
\subset \Coh_\cN(\Loc_\Gv(S))$
and $\QC(\Loc_\Gv(S))
\subset \QC^!_\cN(\Loc_\Gv(S))$. Note as well that 
$\Coh_\cN(\Loc_\Gv(S))$
and $ \QC^!_\cN(\Loc_\Gv(S))$
are respective module categories for  $\Perf(\Loc_\Gv(S))$  and $\QC(\Loc_\Gv(S))$ acting by tensor product.
\end{remark}

Our primary motivation for introducing nilpotent  singular support is the following elementary observation.
Suppose $S$ is a surface with nonempty boundary, and consider the natural projection
$$
\xymatrix{
p:\Loc_\Gv(S, \partial S) \ar[r] & \Loc_\Gv(S)
}
$$
where we forget the $B^\vee$-reduction of $G^\vee$-local systems along $\partial S$. Thus 
it is a base-change of the product over $\pi_0(\partial S)$ of Grothendieck-Springer maps $B^\vee/B^\vee \to G^\vee/G^\vee$. Then 
for any $\cF\in \Coh(\Loc_\Gv(S, \partial S))$, the pushforward
$p_*\cF \in \Coh(\Loc_\Gv(S))$  has nilpotent singular support.
Moreover, every object of $\Coh_\cN(\Loc_\Gv(S))$ arises as such a pushforward, 
and one can extend this to give a full descent description of  $\Coh_\cN(\Loc_\Gv(S))$.

One can argue that the notion of nilpotent singular support is simply a concise way to encode
the idea that we study coherent sheaves coming via parabolic induction.
In the next section, we will describe a global version  of this motivation discovered in the
 de Rham setting.

 %
%

\subsection{Parabolic induction and domain walls}

Any correspondence between classifying spaces of reductive groups
$$
\xymatrix{
B M^\vee &\ar[l]  Z \ar[r] &  B G^\vee
}
$$ 
 induces a correspondence between their character stacks 
$$
\xymatrix{
\Loc_{M^\vee}(S) &\ar[l]   [S,Z] \ar[r] &  \Loc_\Gv(S)
}
$$ 
Depending on technical properties of the maps, one can hope to obtain resulting functors between spectral categories.

From the perspective of topological field theory, 
the correspondence results from the cobordism 
given by the $3$-dimensional cylinder $S\times [0,1]$ labelled with the $M^\vee$ and $\Gv$ theories at the respective ends
$S\times \{0\}$ and $S\times \{1\}$,
and separated by a domain wall or interface  labelled by $Z$ along the separating surface $S\times \{1/2\}$. 

For an important example of this paradigm, fix a parabolic subgroup $P^\vee \subset G^\vee$ with Levi quotient $M^\vee$,
for example a Borel subgroup $B^\vee \subset G^\vee$ with universal Cartan subgroup $H^\vee$.
Arinkin-Gaitsgory~\cite{AG} explain that the choice to work with spectral categories of perfect complexes is not consistent with parabolic induction in the form of geometric Eisenstein series. (They work in the de Rham setting  but the situation is identical in this regard in the Betti setting.) Namely, pulling back and pushing forward 
a perfect complex on $\Loc_{M^\vee}(S)$
 under the natural correspondence
$$\xymatrix{\Loc_{M^\vee}(S) & \ar[l] \Loc_{P^\vee}(S) \ar[r] & \Loc_\Gv(S) 
}$$ 
 leads not to a perfect complex on $\Loc_\Gv(S)$, but at least 
 to a coherent complex with nilpotent singular support.
%
The remarkable theorem of Arinkin-Gaitsgory 
(whose proof in the de Rham setting extends directly to the Betti setting) states that this construction generates all
  coherent sheaves with nilpotent singular support.

\begin{theorem}[\cite{AG}]
The dg category $\QC^!_\cN(\Loc_\Gv(S))$ is generated by the images of
 $\QC(\Loc_{M^\vee}(S))$ under
  the parabolic induction functors 
$$
\xymatrix{
\QC(\Loc_{\Mv}(S))\ar[r] & \QC^!_\cN(\Loc_\Gv(S))
}$$ 
ranging over all
 parabolic subgroups.
\end{theorem}

\begin{remark}The theorem shows the assignment $G^\vee \mapsto \QC^!_\cN(\Loc_\Gv(S))$ is the minimal enlargement of 
the assignment $G^\vee \mapsto  \QC(\Loc_\Gv(S))$ compatible with parabolic induction.
\end{remark}

%
%
%


\subsection{2-sphere and Wilson lines}\label{wilson lines}

Recall that the pushout presentation
$$
\xymatrix{
S^2\simeq D^2 \coprod_{S^1} D^2
}$$
leads to the pullback presentation
$$\Loc_\Gv(S^2)\simeq \{e\}/\Gv \times_{G^\vee/G^\vee} \{e\}/\Gv \simeq \fg^{\vee}[-1]/\Gv$$ 

The underlying classical stack 
$$
\Loc_\Gv(S^2)_{cl} \simeq B\Gv
$$ 
presents the classical Satake category 
$$
\Rep(\Gv) \simeq \Perf(B\Gv) 
$$ 
as the heart of $\Coh(\Loc_\Gv(S^2))$ with respect to the standard $t$-structure.

 The full dg category  $\Coh(\Loc_\Gv(S^2))$ is the  derived Satake category studied in~\cite{BF derived} and also~\cite{AG}. 
Let us  summarize some of its algebraic structures,
in particular its role in spectral Hecke modifications,
 which follow naturally from the perspective of  4-dimensional  topological field theory.
Specifically, if we view $S^2$ as the link of a point in $\RR^4$, then we can view $\Coh(\Loc_\Gv(S^2))$ as line operators,
and seek the algebraic structures they enjoy.

%
%
%
First, by considering  disjoint unions of little 3-disks in a 3-disk, we see that $S^2$ is naturally an $E_3$-algebra in the cobordism category of surfaces. 
Passing to $G^\vee$-local systems,
it follows that
$\Loc_\Gv(S^2)$ is naturally a framed $E_3$-algebra in the correspondence category of derived stacks.
Passing further to coherent sheaves, we obtain an $E_3$-monoidal structure on $\Coh(\Loc_\Gv(S^2))$.

Similarly,
at each  point $s\in S$ of any surface, we see that $S$ carries a natural $E_1$-action  of 
$S^2$ given by little 3-disks in $S\times [0, 1]$ centered along $\{s\} \times (0,1)$.
Passing to $G^\vee$-local systems,
it follows that
$\Loc_\Gv(S)$
carries a natural $E_1$-action  of 
$\Loc_{\Gv}(S^2)$
  in the correspondence category of derived stacks.
  Passing further to coherent sheaves, we conclude that $\Coh(\Loc_\Gv(S))$ is an $E_1$-module over $\Coh(\Loc_\Gv(S^2))$.
  Moreover, there is a natural compatibility between the $E_3$-algebra structure on 
  $\Coh(\Loc_\Gv(S^2))$ and its commuting $E_1$-actions on $\Coh(\Loc_\Gv(S))$   at various points  $s\in S$.

\subsection{3-sphere and local operators}\label{local ops}

The pushout presentation
$$
\xymatrix{
S^3\simeq D^3 \coprod_{S^2} D^3
}
$$
leads to the pullback presentation
$$\Loc_\Gv(S^3)\simeq pt/\Gv \times_{\fg^{\vee}[-1]/\Gv} pt/\Gv \simeq \fg^{\vee}[-2]/\Gv$$ 

 We can view the derived functions  
 $$
 \xymatrix{
 \cO(\Loc_\Gv(S^3)) \simeq  \Sym((\fg^{\vee})^*[2])^\Gv \simeq \Sym(\fh[2])^W
} $$ 
 as local operators
in a 4-dimensional topological field theory.

Similarly as above for line operators,
by considering  disjoint unions of little 4-disks in a 4-disk, we see that $S^3$ is naturally an $E_4$-algebra in the cobordism category of surfaces,
and in fact the endomorphisms of the monoidal unit of $S^2$.
Passing to $G^\vee$-local systems, it follows that
$\Loc_\Gv(S^3)$ is naturally a framed $E_4$-algebra in the correspondence category of derived stacks,
and likewise the endomorphisms of the monoidal unit of $\Loc_\Gv(S^2)$.
Passing further to derived functions, we obtain an $E_4$-monoidal structure on $\cO(\Loc_\Gv(S^3))$,
compatible with its appearance as the endomorphism algebra of the monoidal unit of $ \Coh(\Loc_\Gv(S^2))$.

Going further, at each  point $s\in S$ of any surface,  recall that $ \Coh(\Loc_\Gv(S^2))$ is naturally an $E_1$-module over 
$ \Coh(\Loc_\Gv(S^2))$. Thus we see that $\cO(\Loc_\Gv(S^3))$
acts by endomorphisms of the identity endofunctor of $\Coh(\Loc_\Gv(S))$.
 Moreover, there is a natural compatibility between the $E_4$-algebra structure on 
 $\cO(\Loc_\Gv(S^3))$ and its 
 endomorphisms of the identity endofunctor of $\Coh(\Loc_\Gv(S))$ at various points  $s\in S$.

\begin{remark}

Up to an even grading shift, the local operators $ \cO(\Loc_\Gv(S^3)) $ coincide with the algebra $A$ appearing in the definition
of nilpotent singular support.
In physical language, their  spectrum  provides the  Coulomb branch
$$
 \fg^{\vee}//G^\vee \simeq \fh^\vee//W
$$
of the moduli space of vacua of $\cN=4$ super Yang-Mills theory, with the local operators  the corresponding vacuum expectation values. 
Nilpotent singular support is the natural condition that we sit at the  conformal point $0 \in \fh^\vee//W$ of the moduli of vacua, while considering more general ind-coherent sheaves corresponds to a massive deformation of the theory.

\end{remark}


\subsection{Cylinder and ramification}\label{cylinder}

Similarly as above for line and local operators, there is a  natural concatenation $E_1$-algebra structure on the cylinder
$\cyl$
in the cobordism category of surfaces with boundary. This induces an  $E_1$-algebra structure on 
the Grothendieck-Steinberg stack 
$$
\Loc_\Gv(\cyl, \partial \cyl) \simeq \St_{\Gv} = B^\vee/B^\vee \times_{G^\vee/G^\vee} B^\vee/B^\vee
$$
in the correspondence category of derived stacks. In turn, this induces an
$E_1$-algebra structure on the affine Hecke category  
$$
\cH_{\Gv} = 
\Coh(\St_{\Gv})
$$
compatible with its description~\cite[Theorem 1.4.6(1)]{spectral}  as endofunctors
$$
\xymatrix{
\Coh(\St_{\Gv})
\simeq \End_{\Perf(\Gv/\Gv)}(\Perf(\Bv/\Bv))
}$$

\begin{remark}
There is a subtle difference here  between cocomplete and small categories: passing to ind-categories
is not compatible with taking endofunctors. One could compare the above with the 
parallel statement~\cite{BFN} for cocomplete categories
$$
\xymatrix{
\QC(\St_{\Gv}) \simeq \End_{\QC(\Gv/\Gv)}(\QC(\Bv/\Bv))
}$$
and take note that $\QC(\St_{\Gv})$ is not the ind-completion of $\Coh(\St_{Gv})$.
\end{remark}



Given  any surface with boundary $(S,\partial S)$, and an embedding  $S^1 \subset \partial S$,
 we see that $(S, \partial S)$ carries a natural concatenation $E_1$-action  of 
$\cyl$.
Passing to $G^\vee$-local systems,
it follows that
$\Loc_\Gv(S, \partial S)$
carries a natural $E_1$-action  of 
$\Loc_{\Gv}(\cyl, \partial\cyl)$
  in the correspondence category of derived stacks.
  Passing further to coherent sheaves, we conclude that $\Coh(\Loc_\Gv(S, \partial S))$ is an $E_1$-module over $
  \cH_{\Gv} = \Coh(\Loc_\Gv(\cyl, \partial \cyl ))$.
  Moreover, there are natural compatibilities between these surface operators
  and the previously discussed line operators.
  
  All of the preceding
   equally holds equally well when we restrict to unipotent boundary monodromy and study the  Steinberg stack 
   $$
\Loc^u_\Gv(\cyl, \partial \cyl) \simeq \St^u_\Gv =  \wt\cN^\vee/G^\vee \times_{G^\vee/G^\vee} \wt \cN^\vee/B^\vee
$$
  and the unipotent affine Hecke category
  $$
\cH^u_{\Gv} = 
\Coh( \St^u_\Gv)
$$ 
or to completed unipotent boundary monodromy and study the completed Steinberg stack 
   $$
\Loc^{\hat u}_\Gv(\cyl, \partial \cyl) \simeq \St^{\hat u}_\Gv =  \widehat{\wt\cN^\vee}/G^\vee \times_{G^\vee/G^\vee} \widehat{\wt\cN^\vee}/G^\vee
$$
  and the completed unipotent affine Hecke category
  $$
\cH^{\hat u}_{\Gv} = 
\Coh( \St^u_\Gv)
$$ 

These are the monoidal categoroes 
 appearing in Bezrukavnikov's local geometric Langlands correspondence~\cite{Roma ICM, Roma Hecke},
 see Section~\ref{parabolic bundles} below.
 
 \begin{remark}\label{rem: affine morita}
 The affine Hecke category $\cH_{\Gv}$ is much more complicated and interesting than
 its quasicoherent version $\QC(\St_{\Gv})$. Most notably, thanks to
 Gaitsgory's fundamental 1-affineness theorem~\cite{1affine}, 
  it is possible~\cite{morita} to establish a Morita equivalence 
  $$
  \QC(\St_{\Gv})\module\simeq \QC(G^\vee/G^\vee)\module
  $$
  Thus the noncommutativity of $\QC(\St_{\Gv})$ is essentially trivial (see also Remark~\ref{3Coh}).  
   \end{remark}

\subsubsection{Marked surfaces}\label{spectral ramification}
Given a surface $(S,\partial S)$ with marked boundary, 
we have seen above that  $\Coh(\Loc_\Gv(S,\partial S))$ is naturally
a module for  $\pi_0(\partial S)$-many commuting copies of the affine Hecke category
$\cH_\Gv$. For notational simplicity, we concentrate here on the case of a single marked boundary component $\partial S=S^1$.

On the one hand, given 
any $\cH_\Gv$-module category
$\cM$, 
 we may form a corresponding ramified spectral category
$$\cSpec_\Gv(S,\partial S,\cM)=\Coh(\Loc_\Gv(S,\partial S))\ot_{\cH_\Gv} 
\cM$$
On the other hand, given a reasonable stack $Z\to \Gv/\Gv$, we have a corresponding  stack of ramified local systems
$$\Loc_\Gv(S,\partial S,Z)=\Loc_\Gv(S)\times_{\Loc_\Gv(S^1)} Z$$ which carries a natural singular support condition denoted by $\cN$. 

Now form $\tilde Z = Z  \times_{G^\vee/G^\vee} B^\vee/B^\vee$
and the corresponding $\cH_\Gv$-module $\cM = \Coh(\tilde Z)$.
Then the gluing arguments of~\cite{gluing} extend to the following generality.
%
%

\begin{theorem}\label{gluing ramification}
There is an equivalence of ramified spectral categories 
$$\cSpec_\Gv(S,\partial S,\cM)\simeq \Coh_\cN(\Loc_\Gv(S,\partial S,Z)).$$
\end{theorem}

Thus we can prescribe ramification conditions algebraically or geometrically. Of particular interest is when $Z$
is a moduli of Stokes data for irregular connections on the disk with the map to $\Gv/\Gv$ given by taking monodromy. In this case, the stack $ \Loc_\Gv(S,\partial S,Z)$ is a wild character variety, the Betti version of the moduli of connections with irregular singularities (see~\cite{witten wild} for a discussion of a corresponding wild Geometric Langlands conjecture).


\subsection{$2$-tori and Verlinde loops}
In 3-dimensional Chern-Simons theory, the invariant assigned to
the 2-torus $T^2$ is the Verlinde algebra at level $k$,
which can be realized 
as the fusion ring of characters of  level $k$ integrable representations of the loop group.
The invariant assigned to an arbitrary surface
is the space of WZW conformal blocks.
%
 
 The Verlinde algebra
has a commutative Frobenius algebra structure governing the Verlinde formula for the dimensions of  the 
the spaces of WZW conformal blocks in terms of pants decompositions of surfaces.
It also acts on the spaces of WZW conformal blocks themselves by Verlinde loop operators given by modifications along 
 loops
on surfaces. 

Going to 4-dimensional topological field theory,
an analogous role is played by the  category assigned to $T^2$.
We will outline  this  here and in  Section~\ref{Betti gluing}.  

Recall that $\Loc_{G^\vee}(T^2)$ is the derived commuting stack. 
Along with its more classical realizations, it plays a 
central role in geometric representation theory. Its $K$-theory provides
 a Langlands dual form of the elliptic Hall algebra (see Schiffmann-Vasserot \cite{SV1,SV2,SV3}),
  closely related to Macdonald polynomials and double affine Hecke algebras. 
 It is also directly linked  to Cherednik algebras and the Harish Chandra system or Springer sheaf
 (see Ginzburg \cite{isospectral}).

\subsubsection{Coherent character sheaves}
We refer to the spectral category $\Coh_\cN(\Loc_\Gv(T^2))$ as  the category of {\em coherent character sheaves}.
This name is motivated by the following basic compatibility between the cylinder and $2$-torus.

%
%
%

\begin{theorem}[\cite{spectral}]\label{coherent char sheaves}\label{HH of Hecke} There is an equivalence
$$HH(\cH_\Gv)\simeq \Coh_\cN(\Loc_\Gv(T^2))$$
between the Hochschild homology category of the affine Hecke  category and the category of coherent character sheaves. 
\end{theorem}

\begin{remark}
The theorem is an affine analogue of the result~\cite{character} that the
Hochschild homology of the finite Hecke category of $B$-biequivariant $\D$-modules on $G$
is equivalent to the category of unipotent Lusztig character sheaves on $G$.
\end{remark}

Recall the Hochschild homology category is the home for characters of dualizable module categories. 
Thus any dualizable $\cH_\Gv$-module category 
has a character object in $\Coh_\cN(\Loc_\Gv(T^2))$.
 In particular, for any surface with  boundary $(S,\partial S)$, and embedding
 $S^1\subset\partial S$,
  the spectral category $\Coh_\cN(\Loc_{\Gv}(S,\partial S))$ is a dualizable $\cH_\Gv$-module category, and therefore has a character object in $\Coh_\cN(\Loc_\Gv(T^2))$.
  Its construction is closely related to the global Springer theory of Yun~\cite{Yun} giving actions of Cherednik algebras on cohomology of parabolic Hitchin spaces.

\subsubsection{Verlinde loops}

The pair of pants equips $S^1$ with an $E_2$-algebra structure in the cobordism category
of curves. Taking a product with another $S^1$
  induces the same on the 2-torus $T^2$ 
 in the cobordism category
of surfaces.
Passing to $G^\vee$-local systems,
we obtain an $E_2$-algebra structure on $\Loc_\Gv(T^2)$ in the correspondence category of derived stacks. 

Similarly, at any closed loop $\gamma\subset S$ with trivial normal bundle in a surface, 
we see that $S$ carries a natural $E_1$-action of $T^2$ 
 in the cobordism category
of surfaces.
Passing to $G^\vee$-local systems,
it follows that $\Loc_{\Gv}(S)$ carries a natural  $E_1$-action of $\Loc_{\Gv}(T^2)$
by modifications of $\Gv$-local systems along $\gamma\subset S$.

%
%

One must be careful here when passing to coherent sheaves since not all of the maps in the correspondences
defining the $E_2$-algebra structure on $\Loc_\Gv(T^2)$ are proper.
Let us introduce 
the {\em coherent Verlinde category} $\Coh_{prop/2}(\Loc_\Gv(T^2))$ of  coherent sheaves with proper support relative to the projection $\Loc_{\Gv}(T^2) \to 
\Loc_{\Gv}(S^1) \simeq \Gv/\Gv$ given by restriction to the second loop, where the above algebraic structures are defined with respect
to the first loop.
Then the $E_2$-algebra structure on $\Loc_\Gv(T^2)$ induces the same on
$\Coh_{prop/2}(\Loc_\Gv(T^2))$, and similarly,  $\Coh_{prop/2}(\Loc_\Gv(T^2))$ naturally acts on $\Coh(\Loc_{\Gv}(S))$
 along any closed loop $\gamma\subset S$.

We have the following basic compatibility between the cylinder and $2$-torus.

\begin{theorem}[\cite{spectral}]  There is an $E_2$-monoidal equivalence
$$Z(\cH_{\Gv})\simeq \Coh_{prop/2}(\Loc_\Gv(T^2))
$$
between the center or Hochschild cohomology category of the affine Hecke category 
and the coherent Verlinde category.
\end{theorem}
%


\subsection{Dimensional reduction}\label{dimensional reduction}
Cyclic homology provides an intimate relation between calculus on algebraic varieties and the topology of circle actions. 
Let us briefly recall an instance of this developed in~\cite{conns} and~\cite{TV cyclic} that applies to categories
of sheaves.
From the perspective of topological field theory, it is a form of dimensional reduction,
 or more precisely, of a Nekrasov $\Omega$-background as in~\cite{Nekrasov Witten} (see also~\cite{Witten Atiyah} for a physical discussion of the construction).

Given a category $\cC$ with an $S^1$-action,
consider its equivariant localization or Tate construction  
$$\cC^{S^1,Tate}\simeq \cC^{S^1}\otimes_{\CC[u]} \CC[u,u\inv]
$$ 
where we take $S^1$-invariants and invert the action of the generator of the $S^1$-equivariant cohomology algebra $H^*(BS^1)\simeq\CC[u]$. Since $u$ sits in cohomological degree 2, this results in a 2-periodic or $\ZZ/2$-graded category. One can lift this back to a $\ZZ$-graded category in ``mixed" situations where the $S^1$-action comes from a graded action of the affinization $B\GG_a$ of $S^1$, i.e., lifts to a $B\GG_a\rtimes \GG_m$-action. We refer to this process informally as equivariant localization, see~\cite{reps} for the precise statements (in particular the necessary finiteness condition on the coherent sheaves, Koszul dual to the coherence condition on the $\D$-modules).


Let us consider three surfaces with a natural $S^1$-action: the cylinder $\cyl$, M\"obius strip $\Mob$, and $2$-torus $T^2$, where we rotate one of the loops. The corresponding stacks of $\Gv$-local systems were described in Section~\ref{examples} as the Steinberg stack, the stack of Langlands parameters, and the commuting stack. They inherit an
$S^1$-action, and if we restrict to unipotent monodromy along the rotated loop,
the $S^1$-action comes from a $B\GG_a$-action that in turn lifts to a $B\GG_a\rtimes \GG_m$-action.


 Given an involution $\theta^\vee$ of $\Gv$, with associated $L$-group $G^L$ (see Section~\ref{twisted groups}), let $$\Sigma=\{\sigma\in G^L_{-1}\; :\; \sigma^2=e\}/\Gv$$ be the associated set of involutions of $\Gv$, and for $\sigma\in \Sigma$ let $K^\vee_\sigma=(\Gv)^\sigma$ the corresponding symmetric subgroup.

\begin{theorem}[\cite{reps}] Equivariant localization relates the following categories:
\begin{enumerate}
\item  
$$
\Coh(\St^{\hat u}_\Gv)\leadsto \D_{coh}(\Bv\backslash\Gv/\Bv)
$$
i.e., the (completed unipotent) affine Hecke category reduces to the finite 
equivariant Hecke category (see Section~\ref{parabolic bundles}).

 \item 
 $$\Coh({\mathcal La}_{\Gv}^{\hat u,\theta^\vee})\leadsto
\bigoplus_{\sigma\in\Sigma}\D_{coh}(K^\vee_\sigma\backslash\Gv/\Bv)$$ i.e., coherent sheaves on the (completed unipotent) Langlands parameter space reduce to the categories of geometric Langlands parameters of~\cite{ABV,Soergel} (see Section~\ref{LVS}).

\item 
$$\Coh(\Loc_\Gv^{\hat u}(T^2))\leadsto \D_{coh,\cN}(\Gv/\Gv)
$$ 
i.e., coherent sheaves on the (completed unipotent) commuting stack reduce to adjoint-equivaraint $\D$-modules. 

\end{enumerate}
\end{theorem}

\begin{remark} In the case of the torus, one can check in addition the natural compatibility of singular support for coherent sheaves on loop spaces and for $\D$-modules on the base, so that coherent sheaves with nilpotent singular support (coherent character sheaves) reduce to $\D$-modules with nilpotent characteristic variety, i.e., Lusztig character sheaves (see Section~\ref{char sheaves}).
\end{remark}


\subsection{Interlude: factorization homology}

We briefly digress here to recall a simple method to construct topological field theories from commutative algebras.
Given a symmetric monoidal $\oo$-category, we recall how to build topological field theories
by tensoring $E_\oo$-algebra objects 
 over simplicial sets.
This fits into the formalism of a   construction called factorization homology
that applies more generally to $E_n$-algebra objects as well.
 Following~\cite{BFN}, we will be interested in the symmetric monoidal $\oo$-category of cocomplete
dg categories and in particular dg categories
of quasicoherent sheaves with tensor product.  

%
%

It is a basic fact of commutative algebra that coproducts in the category of commutative algebras
over a base 
are given by relative tensor products over the base (or equivalently, that fiber  products in the category of affine schemes are given by
 spectra of relative tensor products).

 This admits a natural  generalization to homotopical colimits of commutative algebra objects in a symmetric monoidal $\infty$-category $\cC$.  In particular, given a simplicial set $S$, and a commutative algebra $A\in CAlg(\cC)$, one can form the
 {\em factorization homology} 
$$\int_S A= S\ot A\in CAlg(\cC)$$ 
It depends only on the homotopy type of $S$,
and can be calculated as the geometric realization of a natural simplicial commutative algebra given by tensor products of $A$ over the simplices of $S$.
Notably, it satisfies the Mayer-Vietoris gluing property
of a generalized homology theory
$$
\xymatrix{
\int_{S\coprod_T S'} A\simeq \int_S A \ot_{\int_T A} \int_{S'} A
}
$$ 

Now let us restrict from simplicial sets to manifolds.  Then given a commutative algebra $A\in CAlg(\cC)$,
 factorization homology provides a natural $n$-dimensional topological field theory $Z_A$
valued in $\cC$. To a closed $n$-manifold $M^n$, we assign the factorization homology
$$Z_A(M^n) = \int_{M^n} A$$
regarded as a plain object of $\cC$.
To a closed $n-1$-manifold $M^{n-1}$,
we form its factorization homology
regarded as an associative algebra object of $\cC$,
and then assign its module objects
$$ Z_A(M^{n-1}) =  (\int_{M^{n-1}} A)\module$$
Since $A$ is a commutative algebra, $A\module$ is again a commutative algebra under tensor product,
and furthermore, there is a natural equivalence
$$ (\int_{M^{n-1}} A)\module \simeq  \int_{M^{n-1}} (A\module) $$
In this way, we may  continue all the way down to a point,
where we ultimately assign 
$n$-fold iterated module objects
$$ Z_A(pt) =  ((A\module)\cdots)\module$$
The above Mayer-Vietoris gluing property and its natural generalizations translate into the gluing axioms of a topological field theory.

\begin{remark}
More generally, we need not start with a commutative algebra object, but only a framed $E_n$-algebra object,
or a variant  called an $n$-disk algebra which depends on the structure group of the manifold. An important special case is that of associative algebras where factorization homology over the circle recovers Hochschild homology in a form in which
its rotation symmetries are more evident than in the traditional definition.
\end{remark}
%
%
%


\subsection{Quasicoherent spectral categories}

Now let us take $\cC$ to be 
 the symmetric monoidal $\oo$-category of cocomplete
dg categories.

For $X$ a derived stack, its dg category $\QC(X)$ of quasicoherent sheaves equipped with tensor product
 provides a commutative algebra object of $\cC$.
Let us specialize to the classifying stack $X=BG^\vee$
 of a complex reductive group where $\QC(BG^\vee) \simeq \Rep(G^\vee)$. 
 This is a perfect stack in the sense of~\cite{BFN}, and so its factorization homology admits a geometric interpretation.
 Namely, for any homotopy type $S$, we have a natural equivalence
 $$
 \QC(\Loc_\Gv(S))\simeq \int_S \QC(BG^\vee)
 $$ 
 In other words, we may either first  calculate the limit of the mapping stack  $\Loc_\Gv(S)\simeq [S, BG^\vee]$,
 and then pass to quasicoherent sheaves,
 or alternatively first pass to quasicoherent sheaves,  
 and then calculate the colimit of factorization homology over $S$.

Now we may
follow the outline recalled above,  
and  construct a 2-dimensional topological field theory 
 with the assignments:
\begin{enumerate}
\item[$\bullet$] closed surface $M^2$  $\leadsto$ $\QC(\Loc_\Gv(M^2))$
\item[$\bullet$] closed curve $M^1$  $\leadsto \QC(\Loc_\Gv(M^1))\module$
\item[$\bullet$] point $pt$ $\leadsto$ $(\QC(BG^\vee)\module)\module$
\end{enumerate}
Note for a surface $S$ with boundary $\partial S$, we have the assignment
 $$\QC(\Loc_\Gv(S))\in \QC(\Loc_{\Gv}(\partial S))\module$$
 given by the restriction map $\Loc_\Gv(S) \to \Loc_\Gv(\partial S)$

In fact, this 2-dimensional topological field theory 
extends to
a  3-dimensional topological field theory with the additional assignment:
\begin{enumerate}
\item[$\bullet$] closed 3-manifold $M^3$ $\leadsto$ $\cO(\Loc_\Gv(M^3))$
\end{enumerate}
Here for a $3$-manifold $M$ with boundary surface $S$, we have the assignment
 $$p_*\cO(\Loc_\Gv(M))\in \QC(\Loc_{\Gv}(S))$$
 given by the restriction map $p:\Loc_\Gv(M) \to \Loc_\Gv(S)$
More generally, to a 3-dimensional cobordism $M^3:S_1\leadsto S_2$,
 there is the correspondence
$$\xymatrix{ \Loc_\Gv(S_1) &\Loc_\Gv(M)\ar[l]_-{\pi_1}\ar[r]^-{\pi_2}& \Loc_\Gv(S_2)
 }$$
and the theory assigns the functor
$$
\xymatrix{
\pi_{2*}\pi_1^*: \QC(\Loc_\Gv(S_1)) \ar[r] & \QC(\Loc_\Gv(S_2))
}$$

\begin{remark}\label{rem: morita}
Recall the identifications 
$$
\xymatrix{
\Loc_{\Gv}(S^1) \simeq G^\vee/G^\vee
&
\Loc_{\Gv}(\cyl, \partial \cyl) \simeq \St_\Gv \simeq B^\vee/B^\vee \times_{G^\vee/G^\vee}B^\vee/B^\vee 
}
$$ 
and that  there is a Morita equivalence between $\QC(G^\vee/G^\vee)$ 
and the quasicoherent affine Hecke category $\QC(\St_\Gv)$ as discussed in Remark~\ref{rem: affine morita}.

Thus in dimensions 2 and 1, the above topological field theory admits an equivalent formulation
 with the assignments:
\begin{enumerate}
\item[$\bullet$] closed surface $M^2$  $\leadsto$ $\QC(\Loc_\Gv(M^2))$
\item[$\bullet$] closed curve $M^1$  $\leadsto \QC(\Loc_\Gv(M^1 \times [0,1], M^1\times \{0, 1\}))\module$
\end{enumerate}

For a surface with boundary $(S, \partial S)$,
with an identification $S^1 \simeq \partial S$,
 the assignment
 $$\QC(\Loc_\Gv(S))\in \QC(G^\vee/G^\vee)\module$$
corresponds to the assignment
 $$\QC(\Loc_\Gv(S,\partial S))\in \QC(\St_\Gv)\module$$

For another surface with boundary $(S', \partial S')$,
with an identification $S^1 \simeq \partial S'$,
 the gluing axiom
 of the above topological field theory
 $$
\xymatrix{
\QC(\Loc_{\Gv}(S\coprod_{S^1} S')) \simeq \QC(\Loc_{\Gv}(S) \otimes_{\QC(\Gv/\Gv)}
\QC(\Loc_{\Gv}(S') )
}$$
takes the alternative form 
$$
\xymatrix{
\QC(\Loc_{\Gv}(S\coprod_{S^1} S')) \simeq \QC(\Loc_{\Gv}(S, \partial S)) \otimes_{\QC(\St_\Gv)}
\QC(\Loc_{\Gv}(S', \partial S') )
}$$

%
%
\end{remark}


\subsection{Betti spectral gluing}\label{Betti gluing}\label{spectral gluing}
Following the preceding paradigms, we expect the spectral category $\Coh_\cN(\Loc_\Gv(S))$ is the assignment
to the surface $S$ in a 2-dimensional
topological field theory.
Unfortunately, it cannot  be directly described via factorization homology due to the failure of the kind of Morita equivalence appearing in
Remark~\ref{rem: morita}.
Nevertheless, following~\cite{gluing}, it satisfies a Verlinde gluing property which we recall here.

First, recall the affine Hecke category
$$
\cH_{\Gv} = \Coh(\St_\Gv) \simeq \End_{\Perf(\Gv/\Gv)} (\Perf(B^\vee/B^\vee)
$$
We expect the 2-dimensional
topological field theory to assign $\cH_{\Gv}\module$ to the circle.

Next, suppose we have surfaces with boundary $(S, \partial S)$, $(S', \partial S')$,
with identifications $S^1 \simeq \partial S \simeq \partial S'$, then we obtain module categories
$$
\Coh(\Loc_\Gv(S, \partial S)), \hskip.1in \Coh(\Loc_\Gv(S', \partial S')) \in \cH_{\Gv}\module
$$

Finally, we have the following Verlinde gluing property.

\begin{theorem}[\cite{gluing}]
There is a canonical equivalence
$$
\xymatrix{
\Coh_\cN(\Loc_{\Gv}(S\coprod_{S^1} S')) \simeq \Coh(\Loc_{\Gv}(S, \partial S)) \otimes_{\cH_\Gv}
\Coh(\Loc_{\Gv}(S', \partial S')) 
}$$
respecting Hecke symmetries  and Verlinde loop operators.
\end{theorem}

\begin{remark}
The theorem provides motivation for nilpotent singular support:
though it is irrelevant for the punctured surfaces of the right hand side,
 it nevertheless arises for the closed surface of the left hand side.
\end{remark}

A more general version of the theorem holds (see~\cite{gluing}) when the gluing loop is not necessarily separating.
An important special case is when we obtain a 2-torus from a cylinder by gluing its boundary circles as described by Theorem~\ref{HH of Hecke}. 
Applying the general version iteratively to a decomposition of a surface $S$ reduces the calculation of its spectral category 
$\Coh_\cN(\Loc_{\Gv}(S))$ to the fundamental building blocks:
the cylinder,  the disk, the pair of pants, and in the unoriented case, additionally the M\"obius strip.

\begin{remark}\label{3Coh}
To fully construct the  2-dimensional topological field theory, one should find the correct assignment to a point
and invoke the Cobordism Hypothesis.
We conjecture that the 2-category
$\cH_{\Gv}\module$ assigned to the circle is equivalent to ``$2\Coh_\cN(G^\vee/G^\vee)$"
consisting of  smooth categories, proper over 
$G^\vee/G^\vee$,  with nilpotent singular support. 
We further conjecture that the 3-category assigned to a point takes the form
``$3\Coh_\cN(BG^\vee)$" consisting of  categories with smooth diagonal, proper over 
$BG^\vee$, and with nilpotent singular support.
 \end{remark}

\subsection{Quantization}
The character stacks $\Loc_\Gv(S)$ of oriented surfaces carry canonical symplectic structures due to Goldman and Atiyah-Bott, with symplectic pairing on $T_E\Loc_\Gv(S)\simeq H^{*-1}(S,ad(E))$ given by a combination of the cup product and Killing form. (Following~\cite{PTVV} the symplectic form on $\Loc_\Gv(S)=[S,B\Gv]$ derives from the two-shifted symplectic form on $B\Gv$.) 

In~\cite{quantum1,quantum2} the Betti categories $$\QC(\Loc_\Gv(S))\simeq \int_S \Rep(\Gv)$$ were quantized using representation theory of quantum groups. For $q\in\Cx$, let $\Rep_q(\Gv)$ denote the category of algebraic representations of the Drinfeld-Jimbo quantum group associated to $\Gv$. This is a braided tensor category (in fact a balanced tensor category), and thus an $E_2$ (in fact framed $E_2$) algebra in a suitable symmetric monoidal $\infty$-category of categories. In~\cite{quantum1} with Brochier and Jordan we explained how to integrate the quantum group over any oriented surface using factorization homology,
$$S\mapsto \QC_q(\Loc_\Gv(S)):=\int_S \Rep_q(\Gv).$$
These quantum character stacks form a deformation of the topological field theory given by quasi coherent sheaves on character stacks, providing $q$-deformations of the structures we have been discussing.
\begin{remark}
More precisely, we can perform this integration in the setting of dg categories or of abelian categories, and the results of~\cite{SAG} are used to show that the former carry natural t-structures with hearts given by the latter.
\end{remark}
By the characterization of factorization homology, for every point $x\in S$, we have a functor $$\Delta_x: \Rep_q \Gv\to \QC_q(\Loc_\Gv(S))$$ defined by the embedding of any small disc around $x$, and these functors generate the category. For example localizing the trivial representation gives a distinguished object, the quantum structure sheaf $\cO_q\in \QC_q(\Loc_\Gv(S))$.

The categories attached to punctured surfaces are described in~\cite{quantum1} using a $\Rep_q \Gv$ action coming from the inclusion of boundary points: 

\medskip

\begin{theorem}~\cite{quantum1} For any punctured surface $S$ there's an equivalence $\QC_q(\Loc_\Gv(S))\simeq A_q(S)\module_{\Rep_q\Gv}$ with modules for a canonical algebra object in $\Rep_q\Gv$. Moreover a decomposition of the surface provides an explicit presentation of $A_q(S)$. 
\end{theorem}

In particular, the algebra $A_q(S)$ recovers the quantum function algebra $\cO_q(\Gv)$ (with its conjugation action) in the case of an annulus, and the algebra $\D_q(\Gv)$ of quantum differential operators for the punctured torus.

In~\cite{quantum2} the category $\QC_q(\Gv/\Gv)=\int_{S^1}\Rep_q \Gv$ is endowed with two explicit tensor structures, a braided structure identified with the Drinfeld center of $\Rep_q(\Gv)$ quantizing the convolution product on $\QC(\Gv/\Gv)$ and a monoidal structure identified with the trace of $\Rep_q(\Gv)$ and quantizing the pointwise tensor product on $\QC(\Gv/\Gv)$. 
Module categories for the pointwise product on $\QC_q(\Gv/\Gv)$ are then identified with $E_2$-module categories for $\Rep_q(\Gv)$, or with {\em braided} $\Rep_q(\Gv)$-modules. The resulting braided monoidal 
 2-category $\Rep_q(\Gv)\module_{E_2}$ provides the possible markings for surfaces, and the categories associated to surfaces with various markings are then described formally as relative tensor product over $\QC_q(\Gv/\Gv)$, or explicitly as bimodules for $A_q$ of the punctured surface and an algebra describing the marking.
 
\medskip
\begin{theorem}~\cite{quantum2}
The category $\QC_q(\Loc_\Gv(T^2))\simeq \D_q(\Gv/\Gv)$ is identified with adjoint equivariant quantum $\D$-modules on $\Gv$.
 For $\Gv=GL_n$, the category attached to $T^2$ with a ``mirabolic" marking $I_t$ has $End(\cO_q)\simeq {\bf SH}_{q,t}$, the spherical double affine Hecke algebra.
\end{theorem}

\subsubsection{Comparison with the de Rham setting}\label{quantizing Conn}
The de Rham space $\Conn_\Gv(X)$ likewise carries an algebraic symplectic structure, which is analytically identified with that of $\Loc_\Gv(S)$. It is symplectically equivalent to a twisted cotangent bundle of the moduli stack of bundles $Bun_\Gv(S)$, twisted by the determinant line bundle $\fdet$. As a result the category $\QC(\Conn_\Gv(X))$ has a natural deformation quantization, given by modules $\D_{k^\vee}(Bun_\Gv(X))$ over the algebra $\D_{k^\vee}$ of differential operators on $Bun_\Gv(X)$ twisted by $\fdet^{k^\vee}$ (where $k^\vee\in \CC$ is the reciprocal of the quantization parameter). 

The quantum analog of the Riemann-Hilbert correspondence relating Betti and de Rham spaces is provided by the Kazhdan-Lusztig equivalence\cite{KL quantum} between $\Rep_q\Gv$ and the Kazhdan-Lusztig category $KL_{k^\vee}(\Gv)=(\wh{\fg^\vee},\Gv(\cO))\module_{k^\vee}$ of $\Gv(\cO)$-integrable representations of the affine Kac-Moody algebra associated to $\fg^\vee$ at level $k^\vee$. Here the level $k^\vee$ and quantum parameter $q$ are related by $q=\exp(\pi i/(k+h^\vee))$ with $h^\vee$ the dual Coxeter number (with a modification by the ratio of lengths of roots in the non-simply laced case), and $k+h^\vee\notin \QQ_{\geq 0}$.
While $\QC_q(\Loc_\Gv(S))$ is defined as the factorization homology of $\Rep_q \Gv$, the category $\D_{k^\vee}(Bun_\Gv(X))$ is closely related to the corresponding integration (chiral homology) of $KL_{k^\vee}(\Gv)$. In particular for every 
point $x\in X$, Beilinson and Drinfeld constructed localization functors from $$\Delta_x:KL_{k^\vee}(\Gv)\to \D_{k^\vee}(Bun_\Gv(X)),$$ and it was proved by Rozenblyum that these functors generate the latter category on quasi-compact substacks of $Bun_\Gv(X)$.


\section{Automorphic side}

{\em Notation.} For a smooth stack $Z$,  we will write $\Shv(Z)$ for the dg category of complexes of sheaves of vector spaces on the underlying complex analytic stack $Z_{an}$. 
Note we do not require any bounds on the size of such complexes. 
To specify constructible sheaves, we will use the additional notation $\Shv^c(Z) \subset \Shv(Z)$.
For a closed conic subset $\Lambda\subset T^*Z$, we will  write $\Shv_\Lambda(Z)\subset \Shv(Z)$ for the full dg subcategory of those complexes with singular support contained in $\Lambda$. 
%

\subsection{Prelude: character sheaves}\label{char sheaves}
To construct and classify characters of finite groups of Lie type, Lusztig introduced their geometric avatars 
in the form of character sheaves~\cite{character 1} (see also~\cite{laumon character} for a review, and \cite{ginzburg character,MV character,grojnowski} for geometric
approaches to the theory).

Over the complex numbers, character sheaves are perverse sheaves on the adjoint quotient $G/G$ with singular support 
contained in the nilpotent cone $\cN \subset T^*(G/G)$. Recall  that Hamiltonian reduction provides an equivalence
with the moduli of commuting
pairs 
$$T^*(G/G) \simeq \{(g, \xi) \in G, \times \fg^* \, :\,  \on{Ad}_g(\xi) = \xi\} /G
$$
and the nilpotent cone  is the inverse-image of $0\in \fg^*//G$ under the invariant polynomial map 
$$
\xymatrix{
T^*(G/G) \ar[r] & \fg^*//G 
&
(g, \xi) \ar@{|->}[r] & \xi
}
$$

In particular, for a torus $G = T$, the nilpotent cone is the zero-section,
and character sheaves are simply local systems. 
For arbitrary reductive $G$, a natural source of character sheaves is parabolic induction.
Principal series character sheaves arise by starting with such local systems,
and forming their parabolic induction via the correspondence
$$\xymatrix{T/T &\ar[l] B/B \ar[r] & G/G
}$$ 
For example, starting from the trivial local system leads to the Springer sheaf.
A basic observation  is that the characteristic directions of the map $B/B \to G/G$ equal the  nilpotent cone,
and hence parabolic induction preserves nilpotent singular support.
In general, 
all character sheaves arise by parabolic induction from a minimal Levi subgroup, with cuspidal character sheaves
those for which the Levi subgroup is $G$ itself.

\begin{remark}
Over the complex numbers, 
one can equivalently describe character sheaves as $\D$-modules on $G/G$ with nilpotent singular support.
For example, the Springer sheaf and its natural twists form the Harish Chandra system, the eigensystem for the 
bi-invariant differential operators~\cite{HK}. In this case, the condition of nilpotent singular support arises
naturally from the identification of the symbols of the Harish Chandra operators with polynomials on $\fg^*//G$. 
More generally, one can view character sheaves as quantizations of 
fibers of the invariant polynomial map $T^*(G/G)\to \fg^*//G$. In particular, near the identity of $G$, or equivalently on the Lie algebra $\fg$, the Fourier transform identifies character sheaves with coadjoint orbits in $\fg^*$~\cite{mirkovic character}.
With Geometric Langlands in mind, this was interpreted in~\cite{nadler springer} in the language of Fukaya categories, inspired by~\cite{KW}.
\end{remark}


\subsection{Nilpotent sheaves}

Now fix a smooth complex projective curve $X$.

In place of the adjoint quotient $G/G$, we will consider the moduli stack $\Bun_G(X)$ of $G$-bundles on $X$,
and 
in place of the 
invariant polynomial map $T^*(G/G)\to \fg^*//G$,
we will consider the Hitchin integrable system.

%

Set $A_G(X) = H^0(X, (\fg^*//G) \otimes \omega_X)$,
and recall the  role of Higgs bundles 
$$
\xymatrix{
T_E^* \Bun_G(X) \simeq H^0(X, \fg_E^*\otimes \omega_X)
&
E\in \Bun_G(X)}
$$
Then the Hitchin integrable system is the natural map
$$
\xymatrix{
Hitch: T^*\Bun_G(X)\ar[r] &  A_G(X)
& Hitch (E, \phi) = \phi
}
$$ 
induced by $\fg^*\to \fg^*//G$.

The global nilpotent cone, introduced by Laumon~\cite{laumon nilpotent}, is the zero-fiber
$$\cN_{X, G}=Hitch\inv(0)\subset T^*Bun_G(X)
$$
parameterizing $G$-bundles and nilpotent Higgs fields.
It is a conic Lagrangian substack 
with respect to the natural algebraic symplectic structure~\cite{laumon nilpotent,faltings,ginzburg nilpotent, BD} 

With this setup,  natural analogues of character sheaves are Hecke eigensheaves,
geometric avatars of automorphic functions.
Indeed, this perspective in part motivated Laumon's introduction of the global nilpotent cone~\cite{laumon langlands},
and in particular, he conjectured that the global nilpotent cone contained the singular support of all cuspidal Hecke eigensheaves for $GL_n$.

In particular, for a torus $G = T$, the global nilpotent cone is the zero-section,
and nilpotent sheaves are simply local systems. For arbitrary reductive $G$, a natural source of nilpotent sheaves
is parabolic induction. The most basic form of 
geometric Eisenstein series
arise by starting with such local systems, and forming their parabolic induction via the correspondence
$$
\xymatrix{\Bun_T(X) &\ar[l]_-q \Bun_B(X) \ar[r]^-p & \Bun_G(X)
}$$ 
It was explained in~\cite{ginzburg nilpotent} that  the characteristic directions of $p$
equal the global nilpotent cone
$$
\xymatrix{
\cN_{X, G}=\{(E,\phi)\in T^*\Bun_G(X)\; : \; \exists E_B\in p\inv(E) \mbox{ such that } dp^*\phi|_{E_B}=0\}
}$$
and hence such geometric Eisenstein series are nilpotent sheaves.
More generally, one also obtains nilpotent sheaves by applying 
geometric Eisenstein series to nilpotent sheaves for Levi subgroups.

Let us formalize the idea that we should focus on sheaves with nilpotent singular support.
Let us not fix the size requirements of a constructible sheaf, or any parameters for example satisfied by an irreducible Hecke eigensheaf,
but simply require nilpotent singular support.

\begin{definition}
1) We say that $\cF\in \Shv(\Bun_G(X))$ is a {\em nilpotent sheaf} if its singular support lies within 
the global nilpotent cone.

2) Define the large Betti automorphic category 
$$\Shv_{\cN}(Bun_G(X)) \subset \Shv(\Bun_G(X))$$ 
to be the dg category of  nilpotent sheaves,
and 
the small Betti automorphic category 
$$\Shv^w_{\cN}(Bun_G(X)) \subset \Shv_{\cN}(Bun_G(X))$$ 
to be the full dg subcategory of compact objects.  
\end{definition}

\begin{remark}
The superscript $w$ stands for ``wrapped" as appears in the discussion of Section~\ref{nilpotent discussion}.
An object $\cF\in \Shv^w_{\cN}(Bun_G(X))$ is not necessarily constructible, though there is a stratification of $\Bun_G(X)$
such that the restrictions of $\cF$ to the strata are locally constant.
\end{remark}

Here is an initial rough form of the Betti Geometric Langlands conjecture.

\begin{conj}\label{conj: rough}
There is an equivalence
$$\Shv_{\cN}(Bun_G(X)) \simeq \QC^!_\cN(\Loc_\Gv(X))
$$
compatible with Hecke functors (see Section~\ref{Satake} below).
\end{conj}

\begin{remark}\label{integral Betti}
The Betti conjecture admits a natural integral form, relating coherent sheaves on a $\ZZ$-form of the character stack (see Remark~\ref{integral form}) with complexes of sheaves of abelian groups on $Bun_G(X)$, compatibly with actions of the integral form of the Hecke category~\cite{MV satake}.
\end{remark}

%
%
%
%


\subsection{Betti class field theory}
Let $G=T$ be a torus, $\Tv$ the dual torus, and $\Lambda=\Hom(\Cx,T)$ the cocharacter lattice.
The choice of a point $x\in X$ gives rise to an identification
$$\Bun_T(X)\simeq Pic_T(X)^0\times BT\times \Lambda$$
The nilpotent cone $\cN=\{0\}\subset T^*\Bun_T(X)$ is the zero section, so that 
the automorphic category comprises local systems. 

On the other hand, we have an identification 
$$\Loc_{\Tv}(X)= Hom(\pi_1(X),T^\vee) \times B\Tv \times \Spec \Sym \ft[1]$$ 
and the spectral category comprises quasicoherent sheaves $\QC^!_\cN(\Loc_\Tv(X))=\QC(\Loc_\Tv(X))$

We can now use these product decompositions to explicitly match the two categories. The automorphic category is graded by $\Lambda = \pi_0(Bun_T(X))$
matching the grading of the spectral category by $\Lambda = K_0(B\Tv)$.
Next, local systems on $BT$ are given by complete modules over $H^*(BT)=\Sym \ft^*[-2]$, or, via Koszul dually, by quasicoherent modules over $H_*(T)=\Sym \ft[1]$, matching quasicoherent sheaves on the derived factor of $\Loc_{\Tv}(X)$.
Finally, local systems on $Pic_T(X)^0$ are modules for 
$$k[\pi_1(Pic_T(X)^0)]\simeq k[H_1(X)\ot \Lambda],$$ 
as are quasicoherent sheaves on $$Hom(\pi_1(X),T^\vee)= T^\vee \otimes H^1(X)=\Spec( k[H_1(X)\ot \Lambda]).$$



 

\subsection{Hecke modifications}\label{Satake}

Let us recall the action of the Satake category via Hecke modifications
 parallel to the Wilson lines discussed in Section~\ref{wilson lines}. 

Set $\cO = \CC[[t]]$, $\cK = \CC((t))$.
For a point $x\in X$, let $D_x = \Spec \cO_x$ be the disk, where $\cO_x$ is the completed local ring, and   let $D^\times_x = \Spec \cK_x$ be the punctured disk, 
where $\cK_x$ is the fraction field.
For a choice of local coordinate, we obtain isomorphisms $ \cO_x \simeq \cO$,
 $ \cK_x \simeq \cK$.

Set $X_\pm = X$, and introduce  the non-separated curve $X(x) = X_- \coprod_{X\setminus \{x\}} X_+$,
and consider the natural correspondence
$$
\xymatrix{
\Bun_G(X_-) & \ar[l] \Bun_G(X(x)) \ar[r] & \Bun_G(X_+)
}
$$
For a choice of trivialization of a $G$-bundle on $D_x$, 
the corresponding fibers of the maps are isomorphic to the affine Grassmannian $Gr_{G}=G(\cK)/G(\cO)$. It follows the possible kernels for integral transforms of sheaves are $G(\cO)$-equivariant  sheaves
on $Gr_G$.

 Recall the Geometric Satake Theorem~\cite{q-analog,ginzburg satake,BD, MV satake} is an equivalence between $G(\cO)$-equivariant perverse sheaves on $Gr_G$ and finite-dimensional representations of $\Gv$. 
 It admits the following derived enhancement.
 

\begin{theorem}
[\cite{BF derived, AG}] 
There is an equivalence of monoidal dg categories 
$$\Shv_c(G(\cO)\backslash Gr_G)\simeq \Coh(\Loc_\Gv(S^2))
$$
\end{theorem}
%
%

\begin{remark}[twisted local systems]\label{twisted}
An important subtlety in the Geometric Satake Theorem (see~\cite{MV satake,BD} and especially~\cite{reich}) is that the  geometric commutativity constraint on perverse sheaves on $Gr_G$ does not match the algebraic commutativity constraint of representations of $\Gv$.
Rather, it matches that  of  representations of $\Gv$ on super vector spaces in which a canonical central involution $s_\Gv\in Z(\Gv)$ acts by the parity operator  (following the notation of~\cite{FG}). To correct for this, 
one ought to either consider twisted sheaves on  $\Bun_G(X)$, or to  consider $\Gv$-local systems on the $\ZZ/2$-gerbe of spin structures on $X$ with $-1$ acting by $s_\Gv$. This issue can usually be suppressed by choosing a spin structure, and we will only mention it briefly when discussing tame ramification in Section~\ref{parabolic bundles}
and real bundles on unoriented surfaces in Section~\ref{LVS section}.
\end{remark}

Recall for $x\in X$ the action of $\Coh(\Loc_\Gv(S^2))$ on  the spectral category
$\QC^!_\cN(\Loc_\Gv(X))$ by Wilson lines, and
that it only depends on $x\in X$ through its structure as a homotopy point.
  
  Let us make Conjecture~\ref{conj: rough} more explicit by specifying the above Hecke functors
  match the Wilson lines under the Geometric Satake Theorem.
  Thus for Conjecture~\ref{conj: rough} to hold, the following must also be true.
  (We expect a proof to appear in \cite{zhiwei}.)

\begin{conj}\label{conj: sph mods}
The action of the Satake category
$\Shv_c(G(\cO)\backslash Gr_G)$ by Hecke functors at $x\in X$
preserves nilpotent sheaves $\Shv_\cN(\Bun_G(X))$ and is locally constant
as we vary $x\in X$.
\end{conj}


\subsection{Tame ramification}\label{parabolic bundles}
Let $B\subset G$ be a Borel subgroup.

For a point $x \in X$, consider the moduli $\Bun_G(X,  x)$ of $G$-bundles on $X$
equipped with $B$-reductions at $x$.

Set $X_\pm = X$,  and let $  x_\pm\in  X_\pm$ denote $ x \in X$. 
Introduce  the non-separated curve $X( x) = X_- \coprod_{X\setminus \{ x\}} X_+$,
and consider the natural correspondence
$$
\xymatrix{
\Bun_G(X_-,  x_-) & \ar[l] \Bun_G(X( x),  x_- \cup  x_+) \ar[r] & \Bun_G(X_+,  x_+)
}
$$
For a choice of trivialization of a $G$-bundle on $D_x$ with $B$-reduction at $x$, 
the corresponding fibers of the maps are isomorphic to the affine flag manifold $Fl_{G}=G(\cK)/I$
where we write $I\subset G(\cO)$ for the Iwahori subgroup. It follows the possible kernels for integral transforms of sheaves are $I$-equivariant  sheaves
on $Fl_G$.

We have the fundamental Local Langlands Theorem of Bezrukavnikov for the affine Hecke category.
 It categorifies Kazhdan-Lusztig's geometric realization of the affine Hecke algebra as a convolution
algebra in the equivariant K-theory of Steinberg varieties~\cite{KL} (see also~\cite{CG}).

\begin{theorem}[\cite{Roma ICM, Roma Hecke}]\label{roma thm} There is an equivalence of monoidal dg categories 
$$\Shv_c(I\bs Fl_G) \simeq \Coh(\St^u_\Gv)$$
\end{theorem}

\begin{remark}
The theorem is compatible with the Geometric Satake Theorem via 
Gaitsgory's central functor~\cite{dennis central}.
\end{remark}

\begin{remark}
If we replace $I$-equivariant  sheaves
on $Fl_G$ with bimonodromic sheaves, 
the above extends to a natural family of theorems. We conjecture 
one can formulate them all at once in the Betti setting by taking bimonodromic sheaves
with arbitrary monodromy and dropping the unipotent requirement on the Steinberg stack. 
\end{remark}

We have the following natural extension of Conjecture~\ref{conj: rough} to tame ramification.
  
  Recall the  role of Higgs bundles with simple poles
$$
\xymatrix{
T_{(E, E_B)}^* \Bun_G(X, x) \simeq \{\phi\in H^0(X, \fg_E^*\otimes \omega_X(x)) \, :\, res_x(\phi)  \in (\fg/\fb)_{E_B}^* \}
}
$$
for $E\in \Bun_G(X)$ with $B$-reduction $E|_x \simeq E_B$.  Thus we have a global nilpotent cone
$\cN_{G, (X, x)} \subset T^*\Bun_G(X, x)$
parameterizing $G$-bundles on $X$ with $B$-reductions at $x$ and nilpotent Higgs fields.
Let us take  nilpotent sheaves  $\Shv_{\cN}(Bun_G(X, x)) \subset \Shv(Bun_G(X, x))$ to be the full dg subcategory
of sheaves with singular support in the global nilpotent cone.
  
 Introduce the topological surface $S = X \setminus D^\circ$
obtained by removing a small open topological disk $D^\circ\subset X$ around $x\in X$.

  For simplicity, let us assume $G$ is simply connected.
  This allows us to trivialize the twists mentioned in Remark~\ref{twisted}.

\begin{conj}\label{conj: tame rough} Assume $G$ is simply connected.

There is an equivalence
$$\Shv_{\cN}(Bun_G(X, x)) \simeq \QC^!_\cN(\Loc^u_\Gv(S, \partial S))
$$
compatible with the actions of the affine Hecke category.
\end{conj}

\begin{remark}
If we replace $B$-reductions by $N$-reductions,
the above conjecture extends to a natural family of conjectures.
In this setting, the condition of nilpotent  singular support forces  sheaves
to be monodromic for the natural $H$-action.

There is also an evident generalization of the conjecture where we allow more points of ramification. 

\end{remark}


\subsection{Topological strategies}
Recall that the character stack $\Loc_\Gv(X)$, and hence the spectral category
$\QC^!(\Loc_\Gv(X))$ as well,
  is an invariant of  the homotopy type of the curve $X$.
  Thus for Conjecture~\ref{conj: rough} to be true, the following must also hold.
  
\begin{conj}\label{conj: topological}
The Betti category $\Shv_{\cN}(Bun_G(X))$ of nilpotent sheaves
 depends on the curve $X$ only through its underlying topological surface. 
\end{conj}

%

In general, the geometry of the nilpotent cone $\cN_{X, G}$ is sensitive to the algebraic geometry of $X$,
and more specifically, undergoes jumps along Brill-Noether loci in the moduli of curves. Thus the conjecture 
 is far from obvious. (See~\ref{elliptic section} for the nontrivial case of elliptic curves.)

Building upon Conjecture~\ref{conj: topological}, we can hope to establish an automorphic analogue of the spectral gluing
discussed in Section~\ref{spectral gluing}.
 There is a natural geometric mechanism given by the asymptotic degeneration of the loop group~\cite{faltings, solis}.

\begin{conj}\label{conj: gluing}
The Betti category $\Shv_{\cN}(Bun_G(X))$ of nilpotent sheaves
admits a gluing description under the degeneration of $X$ to a nodal curve.
\end{conj}

The conjecture would reduce the challenge of a Betti Langlands correspondence to the building blocks associated to the disc, cylinder and pair of pants (and M\"obius strip in the unoriented case).

\subsection{Building blocks}

We now discuss the status of the Betti conjecture for the basic building blocks, the once, twice and thrice punctured spheres. See Section~\ref{LVS section} for the M\"obius strip.

The moduli stacks $Bun_G(\PP^1,0)$ and $Bun_G(\PP^1,0,\infty)$ associated to the disk and cylinder are very closely related to the quotients $G(\cO)\backslash G(\cK)/I$ and $I\backslash G(\cK)/I$ of the affine flag variety. The resulting categories of sheaves are related by a Radon transform, and therefore are well understood from the work of Arkhipov and Bezrukavnikov~\cite{AB,Roma Hecke} as the standard module and regular bimodule for the affine Hecke category.

\begin{theorem}\cite{AB, Roma Hecke}
The unipotent form of the Betti Conjecture holds for the disk and cylinder, i.e., we have equivalences
$$\Shv_\cN^u(Bun_G(\PP^1,0))\simeq \QC^!_\cN(\Loc_\Gv(D,S^1)),$$
$$\Shv_\cN^u(Bun_G(\PP^1,0,\infty))\simeq \QC^!_\cN(\Loc_{\Gv}(\cyl,\partial \cyl))$$ compatible with affine Hecke module structures. 
\end{theorem}
%

Langlands duality for the thrice punctured sphere appears deep in general. The case of $SL_2$ is distinguished in that the moduli stack of parabolic bundles has a discrete set of isomorphism classes, making it amenable to direct analysis.

\begin{theorem}\cite{zhiwei}
Conjecture~\ref{Betti conjecture} holds for $G=SL_2$ and $X=\PP^1\setminus \{0,1,\infty\}.$
\end{theorem}


\subsection{Elliptic character sheaves}\label{elliptic section}
In this section, we focus on  Betti Langlands for an elliptic curve $E$. 
The moduli stack $Bun_G(E)$ has close connections to the geometric representation theory of $G$. First, if we restrict to the open substack $G_E:=Bun_G(E)^{0,ss}$ of semistable degree $0$ bundles, we obtain the ``elliptic adjoint quotient", which specializes to $G/G$ (when $E$ is a nodal Weierstrass cubic) and of $\fg/G$ (when $E$ is a cuspidal cubic). The global nilpotent cone on $Bun_G(E)$ restricts to the standard nilpotent cone on $\fg/G$ and $G/G$. Thus perverse sheaves in the automorphic category $\Shv_\cN(Bun_G(E))$ restrict to character sheaves on $\fg$ or $G$. 

On the other hand, the entirety of $Bun_G(E)$ is well known to be a model for the geometry of the (somewhat forbidding) adjoint quotient of the loop group $LG/LG$ --  an idea that originates with Looijenga's (unpublished) identification of holomorphic $G$-bundles on the Tate curve $E_q=\Cx/q^{\ZZ}$ with twisted conjugacy classes in loop groups (see~\cite{EFK,BaGi}). Following this logic (see for example~\cite{S2}, which attributes the idea to Ginzburg), we view nilpotent sheaves on $Bun_G(E)$ as a stand-in for nilpotent sheaves on $LG/LG$, i.e., for character sheaves on the loop group. This motivates the following definition.

\begin{definition}
An elliptic character sheaf on $E$ is an object of the automorphic category $\Shv_\cN(Bun_G(E))$. 
\end{definition}


\begin{remark}
In~\cite{penghui}, Looijenga's idea was advanced to provide a complex analytic uniformization of $G_E$ by adjoint quotients of reductive subgroups of the loop group. This leads to a method to prove the topological invariance of elliptic character sheaves.
\end{remark}

%
%
%

Elliptic character sheaves are expected to be avatars for the nascent theory of affine character sheaves~\cite{affine character, BezKazVar}, which are to play the role for depth zero representations of $p$-adic groups that character sheaves do for finite groups of Lie type.
The theory is already extremely rich if we restrict to the subcategory generated through parabolic induction by trivial local systems on $Bun_M(E)$ for Levi subgroups $M$ --- the elliptic Hall category introduced and studied in depth by 
Schiffmann and Vasserot~\cite{S1,S2,SV1,SV2,SV3}. They identify its decategorification, the elliptic Hall algebra (for all $GL_n$ together), with a variant of Cherednik's double affine Hecke algebra and relate it to Macdonald's symmetric functions and $K$-groups of Hilbert schemes of points. In~\cite{SV3} they extend this analysis to prove a decategorified form of the geometric Langlands conjecture in the formal neighborhood of the trivial $\Gv$ local system, for any curve.

In another direction, in~\cite{elliptic} we restrict to $G_E$ and develop the elliptic analog of Springer theory for the Lie algebra $\fg/G$ (cuspidal case) and Lie group $G/G$ (nodal case). Let  $$W_E = (\pi_1 (E)\otimes \pi_1(T)) \rtimes W$$ denote the double affine Weyl group. Note that by the toy model calculation Section~\ref{elliptic toy} we have an equivalence of abelian categories of representations of $W_E$ with $W$-equivariant quasi coherent sheaves on $\Loc_\Tv(E)$, on which the stabilizers act trivially, as well as with local systems on $Bun_T(E)^{0,ss}$. 
Thus the following theorem carries out the Betti Langlands conjecture on the Springer part of the category of elliptic character sheaves:

\medskip
\begin{theorem}\cite{elliptic}
Parabolic induction induces a fully faithful embedding of abelian categories from $\CC[W_E]\module$ to the heart of $\Shv_\cN(G_E)$.
\end{theorem}

\begin{remark}
Fratila~\cite{Fratila} provided an extension of this theory to all components of the semistable locus $Bun_G^{ss}$, replacing $G\supset B\hookrightarrow T$ on each component by a uniquely chosen parabolic $G\supset P\twoheadrightarrow M$ and the Weyl group $W$ by the relative Weyl group of $M$. 
\end{remark}

%
Recall that Theorem 6.9 from~\cite{character} identifies unipotent character sheaves on $G$ with the trace (Hochschild homology) of the finite Hecke category, while Theorem~\ref{coherent char sheaves} likewise identifies [unipotent] coherent character sheaves $\QC^!_\cN(\Loc_\Gv(T^2))$ with the trace of the [unipotent] affine Hecke category $\cH_\Gv$ [respectively $\cH_\Gv^u$].
Bezrukavnikov's Theorem~\ref{roma thm} and its conjectural extension to all monodromies identify the [unipotent] affine Hecke categories for $G$ and $\Gv$. Therefore the Betti conjecture for $E$ (or a unipotent version thereof) reduces to the following purely automorphic statement:

\begin{conj} The Hochschild homology category of the affine Hecke category $\cH_G$ is equivalent to the category of elliptic character sheaves $\Shv_\cN(Bun_G(E))$.
\end{conj}

This conjecture would also provide a natural source of elliptic character sheaves as characters of module categories for the affine Hecke category -- for example, automorphic categories of curves with parabolic structures. The spectral identification of this character as a coherent character sheaf would then provide a geometric analogue of the Arthur-Selberg trace formula (see also~\cite{frenkel ngo}). 

To imitate the proof of the analogous statement for the finite Hecke category~\cite{character}, one needs a version of the horocycle correspondence relating $Bun_G(E)$ and the $I$-orbits on the affine flag variety $I\backslash LG/I$. A natural relation between the two stacks is provided by the degeneration of the Tate curve $E_q \leadsto E_0$ to a nodal elliptic curve and then the passage to its normalization $\PP^1 \simeq \tilde E_0$. This sequence offers an analogue of the horocycle transform for loop groups completely within the setting of finite-dimensional geometry. We expect the corresponding categories to be related via the resulting geometry of degenerations and normalizations.


\subsection{Real Betti Langlands}\label{LVS section}
The Betti conjecture has a natural extension to real curves, and we will focus here on a distinguished case.
Let $(X,\alpha)$ denote a real form of $X$, where $\alpha$ is a complex conjugation of $X$. 

Let $(G,\theta)$ denote a quasi-split real form of $G$.
The involutions $\alpha$ and $\theta$ define a real form of the moduli of $G$-bundles on $X$, namely the stack $\Bun_{G,\theta}(X,\alpha)$ of $G$-bundles on $X$ identified with their pullback under $\alpha$ and $\theta$. 
Given additionally an $\alpha$-invariant finite subset $Y\subset X$,
there is the parabolic stack $Bun_{G,\theta}(X,Y, \alpha)$ of $G$-bundles on $X$ equipped with flags along $Y$
 and compatibly identified with their pullback under $\alpha$ and $\theta$. 

The quasi-split real form $(G,\theta)$ defines an $L$-group $G^L$,  isomorphic to the semidirect product of $\Gv$ with the Galois group $Gal(\CC/\RR)=\ZZ/2$ acting by the algebraic involution $\theta^\vee$ corresponding to the conjugation $\theta$. 
To avoid further discussion of the twisting mentioned in Remark~\cite{twisted}, we will assume the derived group of $G$ is of adjoint type.

\begin{remark} 
A natural extension of the Betti conjecture when $(X,\alpha)$ has no real points relates sheaves on $Bun_{G,\theta}(X,Y, \alpha)$ with coherent sheaves on a character stack of $\alpha$-twisted $Y$-parabolic $G^L$ local systems on $X$. In the case where $X(\RR)$ is nonempty, one should take into account the additional structure to the Hecke modifications along real points studied in~\cite{nadler thesis}, and given a string theory interpretation in~\cite[Section 6]{Gaiotto Witten}.
\end{remark}

Let us focus on the special case of the real curve $(\PP^1,\alpha)$, with the antipodal conjugation $\alpha$,
and tame ramification along $Y =\{0,\infty\}$. Note that the topological quotient $(\PP^1 \setminus\{0, \infty\})/\alpha$
is the open M\"obius strip $\RR\PP^2\setminus \{0 = \infty\}$ (see Section~\ref{examples}).

\begin{conj}[Affine Langlands-Vogan-Soergel Duality] \label{affine LVS}
There is an equivalence $$\Shv_\cN(\Bun_{G,\theta}((\PP^1,0,\infty), \alpha))\simeq \QC^!_\cN(\Loc^u_{G^L,or}(\Mob,S^1))$$ intertwining natural affine Hecke symmetries.
\end{conj}

\begin{theorem}\cite{sl2r} Conjecture~\ref{affine LVS} holds when $G$ is one of the following: a complex group considered as real group, a torus, $SL_2$ or $PGL_2$. 
\end{theorem}

\subsubsection{Langlands-Vogan-Soergel Duality}\label{LVS}
A primary motivation for Conjecture~\ref{affine LVS} is its relevance, developed in~\cite{cyclic,reps}), to the local Langlands program
over the real numbers~\cite{ABV,Soergel}.

Let us apply $S^1$-equivariant localization to the slightly modified monodromic form of the conjecture 
$$\Shv^{mon}_\cN(\Bun_{G,\theta}((\PP^1,0,\infty), \alpha))\simeq \QC^!_\cN(\Loc^{\wh u}_{G^L,or}(\Mob,S^1))$$ 
where we take unipotent monodromic sheaves on the automorphic side and completed unipotent monodromy on the spectral side.

On the automorphic side,
invoking  traditional paradigms of equivariant localization~\cite{GKM} leads to unipotent monodromic sheaves on the fixed point locus, which is a disjoint union of quotient stacks
$$
\xymatrix{
 \coprod_{\iota\in\Theta} G_{\iota}\backslash G/B
&
\Theta=\{\iota \in G\; : \;  \iota \theta(\iota)=1\}/G
}$$ 
where $G_\iota$ denotes the real form of $G$ corresponding to the conjugation defined by $\theta$ and $\iota$.

On the spectral side, invoking results discussed in Section~\ref{dimensional reduction} leads to $\D$-modules
on the disjoint union of quotient stacks
$$
\xymatrix{
\coprod_{\sigma \in \Sigma} 
K^\vee_\sigma\backslash\Gv/\Bv
&
\Sigma=\{\sigma \in G^\vee\; : \;  \sigma \eta(\sigma)=1\}/G^\vee
}
$$
where $K_\sigma^\vee$ denotes the subgroup of $G^\vee$ fixed by the involution defined by $\eta$ and $\sigma$.

Let us write  $D(\Rep^{\hat 0}(G_{\iota}))$ for the dg category of Harish Chandra modules for $G_\iota$ with generalized trivial infinitesimal character, and $D(\Rep^0 (\fg^\vee, K_\sigma^\vee))$ for that of Harish Chandra modules for $(\fg^\vee, K^\vee_\sigma)$ with trivial infinitesimal character.
Applying the respective localization theories of Kashiwara-Schmid \cite{KS,kashiwara real}
and  Beilinson-Bernstein, we arrive at the following two-periodic (unmixed) form of Soergel's conjecture~\cite{Soergel}
(a lift of Vogan's character duality~\cite{Vogan} from $K$-groups to categories):
 
$$
\bigoplus_{\iota\in \Theta}
D(\Rep^{\hat 0}(G_{\iota}))^{per} \longleftrightarrow \bigoplus_{\sigma\in
\Sigma} D(\Rep^0 (\fg^\vee, K_\sigma^\vee))^{per}
$$ 

\begin{remark}
The compatibility of Conjecture~\ref{affine LVS} with affine Hecke symmetries implies a compatibility of the above statement with finite Hecke symmetries.
%
\end{remark}

\subsection{Quantization}\label{quantum section}
In this section we propose a Betti avatar of the Quantum Geometric Langlands Conjecture, introduced in~\cite{Stoyanovsky} following work of Feigin-Frenkel~\cite{FF1,FF2}, see~\cite{dennis whittaker, dennis quantum}.
Recall from Section~\ref{quantizing Conn} that we have determinant line bundles $\fdet\to Bun_G(X), \fdet\to Bun_\Gv(X)$, and that $\Conn_\Gv(X)$ is identified symplectically with the $\fdet$-twisted cotangent bundle to $Bun_\Gv(X)$. Thus by considering $\fdet^{k^\vee}$-twisted $\D$-modules on $Bun_\Gv(X)$ we have a family of categories with a specialization $$\D_{k^\vee}(Bun_\Gv(X))\leadsto\QC(\Conn_\Gv(X))$$
as $k^\vee\to \infty$. On the other hand we can deform $\D(Bun_G(X))$ by considering $\fdet^k$-twisted $\D$-modules, for $k\in \CC$.
The quantum de Rham conjecture proposes an equivalence
$$\xymatrix{\D_k(Bun_G(X)) \ar[r]& \D_{k^\vee}(Bun_\Gv(X))}$$ where the levels are related by $k^\vee=-1/k$. As $k\to 0, k^\vee\to \infty$ we recover the usual de Rham conjecture (up to the subtlety of $\QC$ vs. $\QC^!_\cN$). The twisted geometric Satake equivalence of~\cite{finkelberg lysenko, reich} identifies the (abelian) spherical Hecke categories acting on quantum $\D$-modules as representations of 
a reductive subgroup of $\Gv$, which is trivial for generic $k$. Instead, the conjecture is required to exchange the localization functor $\Delta_x$ from the Kazhdan-Lusztig category of $\wh{\fg}$ at level $k^\vee$ to $\D_{k^\vee}(Bun_\Gv(X))$ discussed in Section~\ref{quantizing Conn} with a Poincar\'e series functor from twisted Whittaker sheaves on the Grassmannian at $x$ to $\D_k(Bun_G(X))$, see~\cite{dennis whittaker}. 

The Betti category $\Shv_\cN(Bun_G(X))$ also admits a natural twisted version, given by monodromic sheaves. Namely for $q\in \Cx$ we consider twisted sheaves on the $\Cx$-gerbe $\fdet^{\log q}$, or concretely sheaves on the total space of the $\Cx$-bundle $\fdet^\times$ which are locally constant with monodromy $q$ along the fibers. This category is the topological counterpart to twisted $\D$-modules, since $q$-twisted constructible sheaves are identified with $\D$-modules on $Bun_G(X)$ twisted by $\fdet^k$ for any choice of $k$ with $\exp(k)=q$. The condition of nilpotent characteristic variety extends to twisted sheaves, so that we can formulate the following:

\begin{conj}[Quantum Geometric Langlands]
For $q\in \Cx$, let $\Shv_{q,\cN}(Bun_G(X))$ denote $q$-monodromic nilpotent sheaves on $\fdet^\times\to Bun_G(X)$. Then there is an equivalence
$$\xymatrix{\Shv_{q,\cN}(Bun_G(X)) \ar[r]& \QC_q(\Loc_\Gv(X))}$$
\end{conj}

\begin{remark}
We expect the conjecture will need a correction related to a notion of singular support for quantum coherent sheaves, which can be detected as in the classical case by the action of local operators $\int_{S^3} \Rep_q\Gv$.
\end{remark}


\begin{thebibliography}{99}

\bibitem[A]{wrapped loops} M. Abouzaid, On the wrapped Fukaya category and based loops.
J. Symplectic Geom. 10 (2012), no. 1, 27-79. 


\bibitem[Ar]{arinkin thesis}  D. Arinkin, On $\lambda$-connections on a curve where $\lambda$ is a formal parameter. Math. Res. Lett. 12 (2005), no. 4, 551-565.

\bibitem[ABV]{ABV}
J. Adams, D. Barbasch, and D. Vogan, The Langlands Classification
and Irreducible Characters for Real Reductive Groups. {\em Progress
in Mathematics} {\bf 104}. Birkhauser, Boston-Basel-Berlin, 1992.

\bibitem[AB]{AB} S. Arkhipov and R. Bezrukavnikov, Perverse sheaves on affine flags and Langlands dual group.
 arXiv:math/0201073.  Israel J. Math. 170 (2009), 135-183.
 
\bibitem[ABG]{ABG} S. Arkhipov, R. Bezrukavnikov and V. Ginzburg, Quantum Groups, 
the loop Grassmannian, and the Springer resolution.
 arXiv:math/0304173. J. Amer. Math. Soc. 17 (2004), no. 3, 595-678. 

\bibitem[AG]{AG} D. Arinkin and D. Gaitsgory, Singular support of coherent sheaves, and the geometric Langlands conjecture, arXiv:1201.6343. Selecta Math. (N.S.) 21 (2015), no. 1, 1-199. 

\bibitem[BG96]{BaGi}
V. Baranovsky and V. Ginzburg, Conjugacy  Classes in Loop Groups and G-Bundles on Elliptic Curves,  Internat. Math. Res. Notices, No 15 (1996)  734-751.


 
\bibitem[BD]{BD} A. Beilinson and V. Drinfeld, Quantization
of Hitchin Hamiltonians and Hecke Eigensheaves. Preprint, available
at math.uchicago.edu/\~{}mitya.

\bibitem[BGS]{BGS} A. Beilinson, V. Ginzburg and W. Soergel, Koszul
duality patterns in representation theory. {\em Jour. Amer. Math.
Soc.}  {\bf 9}  (1996),  no. 2, 473-527.

\bibitem[BBJ1]{quantum1} D. Ben-Zvi, A. Brochier and D. Jordan, Integrating Quantum Groups over Surfaces.  arXiv:1501.04652.

\bibitem[BBJ2]{quantum2} D. Ben-Zvi, A. Brochier and D. Jordan, Quantum Character Varieties and Braided Module Categories. arXiv:1606.04769.

\bibitem[BN1]{cyclic} D. Ben-Zvi and D. Nadler,
Loop Spaces and Langlands Parameters.
arXiv:0706.0322.

\bibitem[BN2]{conns} D. Ben-Zvi and D. Nadler, Loop Spaces and
  Connections.  arXiv:1002.3636. {\em J. of Topology} (2012) 5(2): 377-430.

\bibitem[BN3]{reps} D. Ben-Zvi and D. Nadler, Loop Spaces and
  Representations.  arXiv:1004.5120. {\em Duke Math J.} (2013)
{\bf 162} (9), 1587-1619. 


\bibitem[BN4]{character}
D. Ben-Zvi and D. Nadler, The Character Theory of a Complex Group, arXiv:0904.1247.

\bibitem[BN5]{elliptic} D. Ben-Zvi and D. Nadler, Elliptic Springer Theory. Compositio Mathematica 151 (2015) 1568-1584.
arXiv:1302.7053

\bibitem[BN6]{gluing} D. Ben-Zvi and D. Nadler, Betti Spectral Gluing. e-print arXiv://1602.07379.


\bibitem[BN7]{sl2r} D. Ben-Zvi and D. Nadler, Affine Langlands-Vogan-Soergel Duality for type $A_1$. Preprint, 2016.


\bibitem[BFN]{BFN} D. Ben-Zvi, J. Francis, and D. Nadler, Integral
  transforms and Drinfeld centers in derived algebraic
  geometry. arXiv:0805.0157. J. Amer. Math. Soc. 23 (2010), 909-966. 


\bibitem[BFN2]{morita} D. Ben-Zvi, J. Francis, and D. Nadler,  Morita equivalence for convolution categories: Appendix to \cite{BFN}.
e-print arXiv:1209.0193.


\bibitem[BNP1]{coherent} D. Ben-Zvi, D. Nadler and A. Preygel, Integral Transforms for Coherent Sheaves. arXiv:1312.7164. To appear, Jour. European Math Soc.

\bibitem[BNP2]{spectral} D. Ben-Zvi, D. Nadler and A. Preygel, A spectral incarnation of affine character sheaves. arXiv:1312.7163.


\bibitem[BL]{BL} J. Bernstein and V. Lunts, 
Equivariant sheaves and functors. Lecture Notes in Mathematics, 1578. Springer-Verlag, Berlin, 1994.

\bibitem[Be1]{Roma ICM} R. Bezrukavnikov, Noncommutative counterparts
  of the Springer resolution. arXiv:math.RT/0604445.  International
  Congress of Mathematicians. Vol. II, 1119-1144, Eur. Math. Soc.,
  Z\"urich, 2006.

\bibitem[Be2]{Roma Hecke}  R. Bezrukavnikov, On two geometric realizations of an affine Hecke algebra. arXiv:1209.0403 

\bibitem[BeB]{BezBrav} R. Bezrukavnikov and A. Braverman, Geometric Langlands correspondence for $\D$-modules in prime characteristic: the GL(n) case. Pure Appl. Math. Q. 3 (2007), no. 1, Special Issue: In honor of Robert D. MacPherson. Part 3, 153-179.

\bibitem[BeF]{BF derived} R. Bezrukavnikov and M. Finkelberg, Equivariant Satake category and Kostant-Whittaker reduction.
 arXiv:0707.3799. Mosc. Math. J. 8 (2008), no. 1, 39-72, 183.

\bibitem[BeKV]{BezKazVar} R. Bezrukavnikov, D. Kazhdan and Y. Varshavsky, A categorical approach to the stable center conjecture.
 arXiv:1307.4669.  Ast\'erisque No. 369 (2015), 27-97.
 
 
\bibitem[BY]{BY} R. Bezrukavnikov and Z. Yun, On Koszul duality for Kac-Moody groups. e-print arXiv:1101.1253.
 Represent. Theory 17 (2013), 1-98.

\bibitem[BG]{eisenstein} A. Braverman and D. Gaitsgory, Geometric Eisenstein series. Invent. Math. 150 (2002), no. 2, 287-384. 

\bibitem[CG]{CG} N. Chriss and V. Ginzburg, Representation Theory and
  complex geometry. Birkh\"auser Boston, Inc., Boston,
  MA, 1997. 
  
\bibitem[D]{donagi lecture}  R. Donagi, The geometric Langlands conjecture and the Hitchin system.
Lecture at the US-USSR Symposium in Algebraic Geometry, Univ. of Chicago,
June-July, 1989.
  
\bibitem[DP1]{DP Hitchin} R. Donagi and T. Pantev,  Langlands duality for Hitchin systems. arXiv:math/0604617. Invent. Math. 189 (2012), no. 3, 653-735. 
  
\bibitem[DP2]{DP Hodge} R. Donagi and T. Pantev,  Geometric Langlands and non-abelian Hodge theory. Surveys in differential geometry. Vol. XIII. Geometry, analysis, and algebraic geometry: forty years of the Journal of Differential Geometry, 85-116, Surv. Differ. Geom., 13, Int. Press, Somerville, MA, 2009.

   
\bibitem[DPS]{DPS} R. Donagi, T. Pantev and C. Simpson, In progress.
  
\bibitem[DrG1]{finiteness} V. Drinfeld and D. Gaitsgory, On some finiteness questions for algebraic stacks.  
Geom. Funct. Anal. 23 (2013), no. 1, 149-294.  e-print arXiv:1108.5351.


\bibitem[DrG2]{BunG compact} V. Drinfeld and D. Gaitsgory, Compact generation of the category of $\D$-modules on the stack of $G$-bundles
on a curve. e-pring arXiv:1112.2402.



\bibitem[EFK]{EFK}
P. Etingof, I. Frenkel, and A. Kirillov Jr., Spherical functions on affine Lie groups, Duke Math. J. 80 (1995), 59-90.

\bibitem[Fa1]{faltings} G. Faltings, 
Stable G-bundles and projective connections. 
J. Algebraic Geom. 2 (1993), no. 3, 507-568. 

\bibitem[Fa2]{faltings verlinde} G. Faltings, 
A proof for the Verlinde formula. 
J. Algebraic Geom. 3 (1994), no. 2, 347-374. 

\bibitem[FF1]{FF1} B. Feigin and E. Frenkel, Duality in W-algebras. Internat. Math. Res. Notices 1991, no. 6, 75-82.

\bibitem[FF2]{FF2} B. Feigin and E. Frenkel, Affine Kac-Moody algebras at the critical level and Gelfand-Dikii algebras. Infinite analysis, Part A, B (Kyoto, 1991), 197-215, Adv. Ser. Math. Phys., 16, World Sci. Publ., River Edge, NJ, 1992. 



\bibitem[FL]{finkelberg lysenko} M. Finkelberg and S. Lysenko,
Twisted geometric Satake equivalence. 
J. Inst. Math. Jussieu 9 (2010), no. 4, 719-739. 

 \bibitem[FG]{FG} V. Fock and A. Goncharov, Moduli spaces of local systems and higher Teichm\"uller theory. Publ. Math. Inst. Hautes \'Etudes Sci. No. 103 (2006), 1-211.
 
\bibitem[FST]{FST} S. Fomin, M. Shapiro and D. Thurston,
Cluster algebras and triangulated surfaces. I. Cluster complexes. 
Acta Math. 201 (2008), no. 1, 83-146. 

\bibitem[Fr]{Fratila} D. Fratila,  On the stack of semistable $G$-bundles over an elliptic curve. arXiv:1406.6593.  Math. Ann. 365 (2016), no. 1-2, 401-421.

\bibitem[F1]{frenkel book} E. Frenkel, Langlands correspondence for loop groups. Cambridge Studies in Advanced Mathematics, 103. Cambridge University Press, Cambridge, 2007. xvi+379 pp.

\bibitem[F2]{frenkel bourbaki} E. Frenkel, Gauge theory and Langlands duality. S\'eminaire Bourbaki. Volume 2008/2009. ExposŽs 997-1011. AstŽrisque No. 332 (2010), Exp. No. 1010, ix-x, 369-403.

\bibitem[FG]{frenkel dennis} E. Frenkel and D. Gaitsgory, Local geometric Langlands correspondence and affine Kac-Moody algebras. Algebraic geometry and number theory, 69-260, Progr. Math., 253, Birkh\"auser Boston, Boston, MA, 2006.

\bibitem[FN]{frenkel ngo} E. Frenkel and Ng\^o B. C., Geometrization of trace formulas. Bull. Math. Sci. 1 (2011), no. 1, 129-199.
 
 
\bibitem[FW]{frenkel witten} E. Frenkel and E. Witten, Geometric Endoscopy and Mirror Symmetry.
 arXiv:0710.5939. Communications in Number Theory and Physics  2 (2008), no. 1, 113-283. 

\bibitem[GMN]{GMN} D. Gaiotto, G. Moore and A. Neitzke, Four-dimensional wall-crossing via three-dimensional field theory. Comm. Math. Phys. 299 (2010), no. 1, 163-224.
 

\bibitem[GW]{Gaiotto Witten} D. Gaiotto and E. Witten, S-Duality of Boundary Conditions In N=4 Super Yang-Mills Theory.
 arXiv:0807.3720.  Adv. Theor. Math. Phys. 13 (2009), no. 3, 721-896.
 
 
\bibitem[G1]{dennis central} D. Gaitsgory, Construction of central elements in the affine Hecke algebra via nearby cycles. 
Invent. Math. 144 (2001), no. 2, 253-280. 


\bibitem[G2]{dennis whittaker} D. Gaitsgory,
Twisted Whittaker model and factorizable sheaves. 
Selecta Math. (N.S.) 13 (2008), no. 4, 617-659. 

\bibitem[G3]{DGcat} D. Gaitsgory, Generalities on DG categories. Available at 
http://www.math.harvard.edu/\~{}gaitsgde/GL/


\bibitem[G4]{outline} D. Gaitsgory, Outline of the proof of the geometric Langlands conjecture for GL2. 
Ast\'erisque No. 370 (2015), 1-112.

\bibitem[G5]{1affine} D. Gaitsgory, Sheaves of categories and the notion of 1-affineness. Stacks and categories in geometry, topology, and algebra, 127-225, Contemp. Math., 643, Amer. Math. Soc., Providence, RI, 2015. 

\bibitem[G6]{dennis quantum} D. Gaitsgory, Quantum Langlands Correspondence. arXiv:1601.05279 

\bibitem[Gi1]{ginzburg character} V. Ginzburg, 
Admissible modules on a symmetric space.
Orbites unipotentes et repr\'esentations, III.
Ast\'erisque No. 173-174 (1989), 9-10, 199--255. 

\bibitem[Gi2]{ginzburg satake} V. Ginzburg, Perverse sheaves on a Loop group and Langlands' duality.
 arXiv:alg-geom/9511007.
 
\bibitem[Gi3]{ginzburg nilpotent} V. Ginzburg, The global nilpotent variety is Lagrangian. Duke Math. J. 109 (2001), no. 3, 511-519. 


\bibitem[Gi4]{isospectral} V. Ginzburg, Isospectral commuting variety, the Harish-Chandra D-module, and principal nilpotent pairs. Duke Math. J. 161 (2012), no. 11, 2023--2111.


\bibitem[GS]{GS} A. Goncharov and L. Shen, Geometry of canonical bases and mirror symmetry. Invent. Math. 202 (2015), no. 2, 487-633.

\bibitem[GKM]{GKM}  M. Goresky, R. Kottwitz and R. MacPherson, Equivariant cohomology, Koszul duality, and the localization theorem. Invent. Math. 131 (1998), no. 1, 25-83. 
 
\bibitem[Gr]{grojnowski} I. Grojnowski, Character sheaves on symmetric spaces. MIT thesis 
(1992), available at http://www.dpmms.cam.ac.uk/~groj/papers.html.

\bibitem[GW]{gukov witten} S. Gukov and E. Witten, Gauge Theory, Ramification, And The Geometric Langlands Program. arXiv:hep-th/0612073. Current developments in mathematics, 2006, 35-180, Int. Press, Somerville, MA, 2008. 

 \bibitem[HRV]{HRV} T. Hausel and F. Rodriguez-Villegas,Mixed Hodge polynomials of character varieties. With an appendix by Nicholas M. Katz. Invent. Math. 174 (2008), no. 3, 555-624. 
   
\bibitem[HLRV]{HLRV} T. Hausel, E. Letellier and F. Rodriguez-Villegas, Arithmetic harmonic analysis on character and quiver varieties. Duke Math. J. 160 (2011), no. 2, 323-400. 

\bibitem[dCHM]{dCHM} M. de Cataldo, T. Hausel and L. Migliorini,Topology of Hitchin systems and Hodge theory of character varieties: the case A1. Ann. of Math. (2) 175 (2012), no. 3, 1329-1407. 
   


\bibitem[HK]{HK} R. Hotta and M. Kashiwara, The invariant holonomic system on a semisimple 
Lie algebra.  Invent. Math.  75  (1984),  no. 2, 327--358. 


\bibitem[KW]{KW}
A. Kapustin and E. Witten, ``Electric-Magnetic Duality And The
Geometric Langlands Program,".Commun. Number Theory Phys. 1 (2007), no. 1, 1-236. e-print hep-th/0604151. 

\bibitem[K]{kashiwara} M. Kashiwara, $\D$-modules and microlocal calculus. 
Translated from the 2000 Japanese original by Mutsumi Saito. 
Translations of Mathematical Monographs, 217. 
Iwanami Series in Modern Mathematics. American Mathematical Society, Providence, RI, 2003.


\bibitem[K2]{kashiwara real} M. Kashiwara, Equivariant derived category and representation of real semisimple Lie groups. Representation theory and complex analysis, 137"1¤74, 
Lecture Notes in Math., 1931, Springer, Berlin, 2008. 

\bibitem[KS]{KS} M. Kashiwara and W. Schmid, 
Quasi-equivariant D-modules, equivariant derived category, and representations of reductive Lie groups. Lie theory and geometry, 457-488, 
Progr. Math., 123, Birkh\"auser Boston, Boston, MA, 1994. 


\bibitem[KL1]{KL} D.  Kazhdan and G. Lusztig, Proof of the
  Deligne-Langlands conjecture for Hecke algebras.  Invent. Math. 87
  (1987), no. 1, 153--215.
  
\bibitem[KL2]{KL quantum} D.  Kazhdan and G. Lusztig,   Tensor structures arising from affine Lie algebras.
J. Amer. Math. Soc. 6 (1993), 905-1011 and 7 (1994), 335-453.


\bibitem[La1]{laumon character} G. Laumon, Faisceaux caract\`eres (d'apr\`es Lusztig). 
S\'eminaire Bourbaki, Vol. 1988/89.  Ast\'erisque  No. 177-178  (1989), Exp. No. 709, 231--260.



\bibitem[La2]{laumon langlands} G. Laumon, Correspondance de Langlands g\'eom\'etrique pour les corps de fonctions. Duke Math. J. 54 (1987), no. 2, 309-359.

\bibitem[La3]{laumon nilpotent} G. Laumon, Un analogue global du c\^one nilpotent. Duke Math. J. 57 (1988), no. 2, 647-671. 

\bibitem[LN]{penghui} P. Li and D. Nadler,  Uniformization of semistable bundles on elliptic curves. arXiv:1510.08762. 


\bibitem[L1]{topos} J. Lurie, Higher topos theory.
arXiv:math.CT/0608040.  Annals of Mathematics Studies, 170. Princeton University Press, Princeton, NJ, 2009.

\bibitem[L2]{HA} J. Lurie, Higher Algebra. Available at
  http://www.math.harvard.edu/\~{}lurie/


\bibitem[L3]{jacob TFT} J. Lurie, On the classification of topological
  field theories.  Available at http://www.math.harvard.edu/\~{}lurie/
  Current developments in mathematics, 2008, 129--280, Int. Press,
  Somerville, MA, 2009.
  
  
\bibitem[L4]{SAG} J. Lurie, Spectral Algebraic Geometry. Available at http://www.math.harvard.edu/\~{}lurie/
  
\bibitem[Lu1]{q-analog} G. Lusztig, 
Singularities, character formulas, and a q-analog of weight multiplicities. Analysis and topology on singular spaces, II, III (Luminy, 1981), 208-229, 
Ast\'erisque, 101-102, Soc. Math. France, Paris, 1983. 
  
\bibitem[Lu2]{character 1} G. Lusztig, Character sheaves I. Adv. Math 56 (1985) no. 3, 193-237.

\bibitem[Lu3]{affine character} G. Lusztig, Unipotent almost characters of simple p-adic groups. arXiv:1212.6540.  Ast\'erisque No. 370 (2015), 243-267. 


\bibitem[M]{mirkovic character} I. Mirkovi\'c, Character sheaves on reductive Lie algebras. Mosc. Math. J. 4 (2004), no. 4, 897-910, 981. 

\bibitem[MV1]{MV character} I. Mirkovi\'c and K. Vilonen, 
Characteristic varieties of character sheaves.  Invent. Math.  93  (1988),  no. 2, 405--418.

\bibitem[MV2]{MV satake} I. Mirkovic and K. Vilonen, Geometric Langlands duality and representations of algebraic groups over commutative rings. Ann. of Math. (2) 166 (2007), no. 1, 95-143.
  
  \bibitem[N1]{nadler thesis} D. Nadler, Perverse sheaves on real loop Grassmannians. Invent. Math. 159 (2005), no. 1, 1-73. arXiv:math/0202150

\bibitem[N2]{nadler springer} D. Nadler, Springer theory via the Hitchin fibration. arXiv:0806.4566. Compos. Math. 147 (2011), no. 5, 1635-1670. 



\bibitem[N3]{Nwms} D. Nadler, Wrapped microlocal sheaves on pairs of pants.
 arXiv:1604.00114.

\bibitem[NY]{zhiwei} D. Nadler and Z. Yun, Preprint, 2016.


\bibitem[NW]{Nekrasov Witten} N. Nekrasov and E. Witten, The Omega Deformation, Branes, Integrability, and Liouville Theory.
 arXiv:1002.0888. J. High Energy Phys. 2010, no. 9, 092, i, 82 pp. 

\bibitem[PTVV]{PTVV} T. Pantev, B. To\"en, M. Vaqui\'e and G. Vezzosi, Shifted symplectic structures. Publ. Math. Inst. Hautes \'Etudes Sci. 117 (2013), 271-328. 


\bibitem[Re]{reich} R. Reich, Twisted geometric Satake equivalence via gerbes on the factorizable Grassmannian. Represent. Theory 16 (2012), 345-449.

\bibitem[S1]{S1} O. Schiffmann, Spherical Hall algebras of curves and Harder-Narasimhan stratas. J. Korean Math. Soc. 48 (2011), no. 5, 953-967. 
 

\bibitem[S2]{S2} O. Schiffmann, On the Hall algebra of an elliptic curve, II. Duke Math. J. 161 (2012), no. 9, 1711-1750. 


\bibitem[SV1]{SV1} O. Schiffmann and E.Vasserot, The elliptic Hall algebra and the equivariant K-theory of the Hilbert scheme of $\mathbb{A}^2$.
arXiv:0905.2555. Duke Math. J. 162 (2013), no. 2, 279-366. 
 
\bibitem[SV2]{SV2} O. Schiffmann and E.Vasserot,The elliptic Hall algebra, Cherednik Hecke algebras and Macdonald polynomials. arXiv:0802.4001. Compos. Math. 147 (2011), no. 1, 188-234.  




\bibitem[SV3]{SV3} O. Schiffmann and E.Vasserot, Hall algebras of curves, commuting varieties and Langlands duality. arXiv:1009.0678. Math. Ann. 353 (2012), no. 4, 1399-1451. 



\bibitem[STWZ]{STWZ} V. Shende, D. Treumann, H. Williams and E. Zaslow, Cluster varieties from Legendrian knots.
 arXiv:1512.08942 

\bibitem[Si1]{simpson} C. Simpson, The Hodge filtration on nonabelian cohomology. Algebraic geometryÑSanta Cruz 1995, 217-281, Proc. Sympos. Pure Math., 62, Part 2, Amer. Math. Soc., Providence, RI, 1997. 


\bibitem[Si2]{simpson opers} C. Simpson,  Iterated destabilizing modifications for vector bundles with connection. 
Vector bundles and complex geometry, 183-206, Contemp. Math., 522, Amer. Math. Soc., Providence, RI, 2010.

\bibitem[S]{Soergel} W. Soergel, Langlands' philosophy and Koszul duality,
      Algebra--representation theory (Constanta, 2000), {\em NATO Sci. Ser.
      II Math. Phys. Chem.}, vol. {\bf 28}, Kluwer Acad. Publ., Dordrecht,
      2001, pp. 379-414.

\bibitem[St]{Stoyanovsky} A. Stoyanovsky, Quantum Langlands duality and conformal field theory. arXiv:math/0610974.

\bibitem[TV]{TV cyclic} B. To\"en and G. Vezzosi, Caract\`eres de Chern, traces \'equivariantes et g\'eom\'etrie alg\'ebrique d\'eriv\'ee. Selecta Math. (N.S.) 21 (2015), no. 2, 449-554.


\bibitem[V]{Vogan} D. Vogan,
Irreducible characters of semisimple Lie groups. IV.
Character-multiplicity duality.  {\em Duke Math. J.}  {\bf
  49}  (1982), no. 4, 943--1073.


\bibitem[W1]{witten wild} E. Witten,
Gauge Theory And Wild Ramification.
 arXiv:0710.0631. Anal. Appl. (Singap.) 6 (2008), no. 4, 429-501. 

\bibitem[W2]{Witten Atiyah} E. Witten, Geometric Langlands and the
  equations of Nahm and Bogomolny.  e-print
  arXiv:0905.4795. Proc. Roy. Soc. Edinburgh Sect. A 140 (2010),
  no. 4, 857-895.

\bibitem[W3]{Witten revisit} E. Witten, More On Gauge Theory And Geometric Langlands. arXiv:1506.04293.

\bibitem[Y]{Yun} Z. Yun, Global Springer Theory. arXiv:0810.2146. Adv. Math. 228 (2011), no. 1, 266-328. 

\end{thebibliography}
\end{document}